\documentclass[11pt]{article}
\usepackage[T2A]{fontenc}
\usepackage[cp1251]{inputenc}
\usepackage{amsmath}
\usepackage{amssymb}
\def\a{\mathbf{a}}
\def\m{\mathbf{m}}
\def\f{\mathbf{f}}

\def\b{\mathbf{b}}
\def\N{\mathbb{N}}
\def\Z{\mathbb{Z}}

\def\C{\mathbb{C}}
\def\O{\mathcal{O}}
\def\aaaa{a}
\def\bbbb{b}
\def\cccc{c}
\def\dddd{d}
\def\eeee{e}
\def\etta{\boldsymbol{\eta}}
\def\zetta{\boldsymbol{\zeta}}
\def\llambda{\boldsymbol{\lambda}}
\def\ggamma{\boldsymbol{\gamma}}
\newtheorem{lemma}{\hspace*{\parindent}Lemma}
\newtheorem{theorem}{\hspace*{\parindent}Theorem}
\newtheorem{corollary}{\hspace*{\parindent}Corollary}
\title{Degenerate Miller-Paris transformations}
\author{D.B.\:Karp$^{\rm a,b}$\footnote{Corresponding author. E-mail: D. Karp -- \emph{dimkrp@gmail.com},
E.\:Prilepkina --  \emph{pril-elena@yandex.ru}}~~and
E.G.\:Prilepkina$^{\rm a,b}$
\\[10pt]
\\
\small{\textit{$\phantom{1}^a$Far Eastern Federal University, 8
Sukhanova street, Vladivostok, 690950, Russia}}
\\
\small{\textit{$\phantom{1}^b$Institute of Applied Mathematics,
FEBRAS, 7 Radio Street, Vladivostok,  690041, Russia}}}
\date{}

\begin{document}
\maketitle

\begin{center}
\parbox{12cm}{
\small\textbf{Abstract.}
Important new transformations for the generalized hypergeometric functions with integral parameter differences have been discovered some years ago by Miller and Paris
and studied in detail in a series of papers by a number of authors. These transformations fail if the free bottom parameter is greater than a free top
parameter by a small positive integer. In this paper we fill this gap in the theory of Miller-Paris transformations by computing the limit cases of these
transformations in such previously prohibited situations.  This leads to a number of new transformation and summation formulas including extensions of Karlsson-Minton theorem.}
\end{center}

\bigskip


\bigskip


\bigskip

\section{Introduction and preliminaries}

Transformation, reduction and summation formulas for the (generalized) hypergeometric functions is a vast subject with rich history dating back to Euler.
Among  important applications of such formulas (let alone hypergeometric functions in general) are quantum physics  \cite{RJRJR}, non-equilibrium statistical
physics \cite{KMT2018} and many other fields \cite{Seaborn}. The main developments up to the end of 20th century can be found, for instance, in the books \cite{AAR,BealsWong,LukeBook}. Over the last decade important new discoveries have been made in this area due, mainly, to Miller and Paris (with important contributions by Rathie, Kim and others), which culminated in the seminal paper \cite{MP2013}.  See also related developments in \cite{KRP2015,Miller2005,MP2011,MP2012,MP2012MathComm}.
These new discoveries deal with the generalized hypergeometric functions with integral parameters differences,
i.e. containing some pairs of top and bottom parameters differing by nonnegative integers.  In the subsequent work \cite{KRP2014} the authors applied the beta
integral method developed and fully automated some years earlier by Krattenthaler and Rao \cite{KrRao2003} to generate Thomae-like relations connecting the generalized
hypergeometric functions of unit argument.

To give a brief account of some of the results from \cite{MP2013} let us introduce the standard notation.  The symbol ${_{p}F_q}(\a;\b;z)$  stands for the generalized hypergeometric function (see \cite[Section~2.1]{AAR}, \cite[Section~5.1]{LukeBook}, \cite[Sections 16.2-16.12]{NIST} or \cite[Chapter~12]{BealsWong}).  When evaluated at $z=1$ the argument of the generalized hypergeometric function will be routinely omitted: ${_{p}F_q}(\a;\b):={_{p}F_q}(\a;\b;1)$. The size of a vector $\a$ is typically obvious from the subscript of the corresponding hypergeometric function.  In the sequel we will also use the shorthand notation for  products and sums:
$$
\Gamma(\a)=\Gamma(a_1)\Gamma(a_2)\cdots\Gamma(a_p),~~(\a)_n=(a_1)_n(a_2)_n\cdots(a_p)_n,~~\a+\mu=(a_1+\mu,a_2+\mu,\dots,a_p+\mu);
$$
in particular, $(\a)=(\a)_1=a_1\cdots{a_p}$; inequalities like $\Re(\a)>0$ and properties like $-\a\notin\N_0$ will be understood element-wise.  The symbol $\a_{[k]}$ stands for the vector $\a$ with omitted $k$-th element.

Let $\m=(m_1,\ldots,m_r)\in\N^r$, $m=m_1+m_2+\ldots+m_r,$ $\f=(f_1,\ldots,f_r)\in\C^r$, where $\N$ and $\C$ stand for the natural and complex numbers, respectively.
The main Miller-Paris transformations are given in \cite[Theorems~3~and~4]{MP2013}. The first generalizes the first Euler (or Euler-Pfaff) transformation for the Gauss hypergeometric
function as follows:
\begin{equation}\label{eq:KRPTh1-1}
{}_{r+2}F_{r+1}\left.\!\!\!\left(\begin{matrix}\aaaa, \bbbb, \f+\m\\\cccc,\f\end{matrix}\right\vert x\right)
=(1-x)^{-\aaaa}{}_{m+2}F_{m+1}\left.\!\!\!\left(\begin{matrix}\aaaa, \cccc-\bbbb-m, \zetta+1\\\cccc, \zetta\end{matrix}\right\vert\frac{x}{x-1}\right).
\end{equation}
It is true for $\bbbb\ne{f_j}$, $j=1,\ldots,r$,  $(\cccc-\bbbb-m)_{m}\ne0$ and $x\in\C\!\setminus\![1,\infty)$, if the generalized hypergeometric function is understood
as analytic continuation for the values of the argument, where the series definition diverges. The second transformation is the confluent version of the first
and generalizes the Kummer transformation for the confluent hypergeometric function. It is given by
\begin{equation}\label{eq:Kummer-type}
{}_{r+1}F_{r+1}\left.\!\!\!\left(\begin{matrix}\bbbb, \f+\m\\\cccc,\f\end{matrix}\right\vert x\right)
=e^{x}{}_{m+1}F_{m+1}\left.\!\!\!\left(\begin{matrix}\cccc-\bbbb-m, \zetta+1\\\cccc, \zetta\end{matrix}\right\vert-x\right)
\end{equation}
and is true for $\bbbb\ne{f_j}$, $j=1,\ldots,r$, $(\cccc-\bbbb-m)_{m}\ne0$ and all complex $x$. The third transformation generalizes the second Euler transformation for the Gauss
hypergeometric function as follows:
\begin{equation}\label{eq:KRPTh1-2}
{}_{r+2}F_{r+1}\left.\!\!\!\left(\begin{matrix}\aaaa, \bbbb, \f+\m\\\cccc,\f\end{matrix}\right\vert x\right)
=(1-x)^{\cccc-\aaaa-\bbbb-m}{}_{m+2}F_{m+1}\left.\!\!\!\left(\begin{matrix}\cccc-\aaaa-m, \cccc-\bbbb-m, \etta+1\\\cccc, \etta\end{matrix}\right\vert x\right)
\end{equation}
and is true for $(\cccc-\aaaa-m)_{m}\ne0$, $(\cccc-\bbbb-m)_{m}\ne0$, $(1+\aaaa+\bbbb-\cccc)_m\ne0$ and $x\in\C\!\setminus\![1,\infty)$.
Here  $\zetta=\zetta(\bbbb,\cccc,\f)=(\zeta_1,\ldots,\zeta_m)$ are the roots of the characteristic polynomial
\begin{equation}\label{eq:Qm}
Q_m(t)=Q_m(\bbbb,\cccc,\f|t)=\frac{1}{(\cccc-\bbbb-m)_{m}}\sum\limits_{k=0}^{m}(\bbbb)_kC_{k,r}(t)_{k}(\cccc-\bbbb-m-t)_{m-k}
\end{equation}
with $C_{0,r}=1$, $C_{m,r}=1/(\f)_{\m}=[(f_1)_{m_1}\cdots(f_r)_{m_r}]^{-1}$, and
\begin{equation}\label{eq:Ckr}
C_{k,r}(\f)=\frac{1}{(\f)_{\m}}\sum\limits_{j=k}^m\sigma_j\mathbf{S}_j^{(k)}=\frac{(-1)^k}{k!}{}_{r+1}F_{r}\!\left(\begin{matrix}-k,\f+\m\\\f\end{matrix}\right),
\end{equation}
where $\mathbf{S}_j^{(k)}$ are the Stirling's numbers of the second kind. The numbers $\sigma_j$ ($0\leq{j}\leq{m}$) are defined via the generating function
$$
(f_1+x)_{m_1}\ldots(f_r+x)_{m_r}=\sum\limits_{j=0}^m\sigma_jx^j.
$$
Further, $\etta=\etta(\bbbb,\cccc,\f)=(\eta_1,\ldots,\eta_m)$ are the roots of the second characteristic polynomial
\begin{equation}\label{eq:Qmhat}
\hat{Q}_m(t)=\hat{Q}_m(\bbbb,\cccc,\f|t)=\sum\limits_{k=0}^{m}\frac{(-1)^kC_{k,r}(\aaaa)_k(\bbbb)_k(t)_k}{(\cccc-\aaaa-m)_k(\cccc-\bbbb-m)_k}
{}_{3}F_{2}\!\left(\begin{matrix}-m+k,t+k,\cccc-\aaaa-\bbbb-m\\\cccc-\aaaa-m+k,\cccc-\bbbb-m+k\end{matrix}\right).
\end{equation}
All three transformations (\ref{eq:KRPTh1-1})-(\ref{eq:KRPTh1-2}) fail when $(\cccc-\bbbb-m)_{m}=0$, that is when $\cccc-\bbbb\in\{1,\ldots,m\}$.
Note, however, that for these values of parameters the left hand sides of (\ref{eq:KRPTh1-1})-(\ref{eq:KRPTh1-2}) remain  well-defined.  Hence, it is natural to ask: what happens to the right hand sides of (\ref{eq:KRPTh1-1})-(\ref{eq:KRPTh1-2}) when $(\cccc-\bbbb-m)_{m}=0$?  The main purpose of this paper is to give a complete answer to this question.

Before we turn to the general situation in the next section, let us make the following  simple observation. The case $r=m=1$, $\cccc=\bbbb+1$
is easy to treat using the straightforward limit
\begin{equation}\label{eq:limit-m1}
\lim\limits_{\varepsilon\to0}{}_{p+1}F_{p}\left.\!\!\left(\begin{matrix}\varepsilon,\a\\\alpha\varepsilon,\b\end{matrix}\right|x\right)
=1-\frac{1}{\alpha}+\frac{1}{\alpha}{}_{p}F_{p-1}\left.\!\!\left(\begin{matrix}\a\\\b\end{matrix}\right|x\right).
\end{equation}
When applied to the $r=m=1$ case of (\ref{eq:KRPTh1-2}) (given explicitly in \cite[(4.11a)]{MP2011}) this approach yields
$$
(1-x)^{\bbbb}{}_{3}F_{2}\left.\!\!\left(\begin{matrix}\aaaa,\bbbb,f+1\\\aaaa+1,f\end{matrix}\right|x\right)=
1-\frac{\bbbb(\aaaa-f)}{f(\aaaa-\bbbb)}+\frac{\bbbb(\aaaa-f)}{f(\aaaa-\bbbb)}{}_{2}F_{1}\left.\!\!\left(\begin{matrix}1,\aaaa-\bbbb\\\aaaa+1\end{matrix}\right|x\right).
$$
Note that ${}_2F_1$ on the right hand side is incomplete beta function \cite[(8.17.8)]{NIST}.
Similarly, the $r=m=1$ case of (\ref{eq:KRPTh1-1}) (given explicitly in \cite[(1.3a)]{MP2011}) reduces to
$$
(1-x)^{\bbbb}{}_{3}F_{2}\left.\!\!\left(\begin{matrix}\aaaa,\bbbb,f+1\\\aaaa+1,f\end{matrix}\right|x\right)
=1-\frac{f-\aaaa}{f}+\frac{f-\aaaa}{f}{}_{2}F_{1}\left.\!\!\left(\begin{matrix}1,\bbbb\\\aaaa+1\end{matrix}\right|\frac{x}{x-1}\right).
$$
Certainly, the two above formulas can be obtained from one another by Euler-Pfaff's transformation \cite[(2.2.6)]{AAR} and the appropriate contiguous relation.
Although it's hard to image that these simple identities could be new we were unable to immediately locate them in the literature.

The remaining part of this  paper is organized as follows.  In the next section we give a complete treatment of the general case $\cccc-\bbbb\in\{1,\ldots,m\}$.
In section~3 we explore the most important extreme cases $\cccc-\bbbb=1,2,m$ when many formulas simplify significantly. The first case $\cccc-\bbbb=1$ can be viewed
as a generalization of Karlsson's summation theorem when argument unity is replaced by an arbitrary argument.
In section~4 we apply asymptotic approximations and the beta integral method to obtain some summation and transformation theorems
the generalized hypergeometric functions at unit argument.  In particular, we give another
extension of Karlsson's summation theorem and some examples of summation theorems for ${}_{4}F_{3}(1)$ with one positive and two negative integral parameters differences.

\section{The general case}

The key to understanding of what happens with (\ref{eq:KRPTh1-1}) and (\ref{eq:KRPTh1-2})
when $\cccc-\bbbb\in\{1,\ldots,m\}$ is provided by the exploration of the limit behavior of the characteristic polynomials $Q_m$ and $\hat{Q}_m$. This is done in the following two lemmas.

\begin{lemma}\label{lm:Qlimit}
Write $Q_m^\varepsilon(t)$ for the polynomial \emph{(\ref{eq:Qm})} with $\cccc-\bbbb-m=-q+\varepsilon$, where $q\in\{0,1,\ldots,m-1\}$, $\varepsilon>0$
and let $\zeta_1(\varepsilon),\ldots,\zeta_m(\varepsilon)$ denote its zeros. Suppose $f_j-\bbbb\notin\{m-q-m_j,m-q-m_j+1,\ldots,0\}$ for indices $0\le{j}\le{r}$ such that $m-q-m_j\le0$. Then the next claims are true as $\varepsilon\to0$\emph{:}

\smallskip

\emph{(1)} the first \emph{(}after possible renumbering\emph{)} $q+1$ zeros tend to consecutive non-positive integers\emph{:} $\zeta_1\to0$, $\zeta_2\to-1$, $\ldots$, $\zeta_{q+1}\to-q$\emph{;}

\smallskip

\emph{(2)} the remaining $m-q-1$ zeros \emph{(}no such zeros remain if $q=m-1$\emph{)}  tend to the zeros of the reduced polynomial
\begin{equation}\label{eq:R}
R_{m-q-1}(t)=\sum\limits_{k=0}^{m}(\bbbb)_{k}C_{k,r}(-1)^{k}(1-t-k)_{m-q-1}
\end{equation}
of degree $m-q-1$\emph{;}

\smallskip

\emph{(3)} the limit relation
\begin{equation}\label{eq:A00Q}
\lim\limits_{\varepsilon\to0}\frac{\varepsilon}{\zeta_1(\varepsilon)}=\frac{R_{m-q-1}(0)}{(m-q-1)!}
\end{equation}
holds\emph{;}

\smallskip

\emph{(4)}  the values $Q_m^\varepsilon(-l)$ at first $q+1$ non-positive integers converge \emph{(}as $\varepsilon\to0$\emph{)} to the following expressions\emph{:}
\begin{equation}\label{eq:Qm0negative}
\lim\limits_{\varepsilon\to0}Q^{\varepsilon}_m(-l)=Q_m^0(-l)=\frac{1}{(-q)_l}\sum\limits_{k=0}^{l}(\bbbb)_kC_{k,r}(-l)_{k}(m-q)_{l-k},
\end{equation}
where $l=0,1,\ldots,q$.  In particular, $Q_m^0(0)=1$.
\end{lemma}
\textbf{Remark.}  In our concurrent paper \cite[Theorem~4]{KPNew} we show that a simpler formula holds for the polynomial $R_{m-q-1}$, namely:
\begin{equation}\label{eq:Rnew}
R_{m-q-1}(t)=\sum\limits_{k=0}^{m-q-1}\frac{(\f-\bbbb-k)_{\m}(\bbbb)_{k}}{(\f)_{\m}k!}(q+1-m)_{k}(\bbbb+k+1-t)_{m-q-1-k}.
\end{equation}
This formula shows that if $m-q-m_j\le0$ for some $0\le{j}\le{r}$ and $f_j-\bbbb\in\{m-q-m_j,m-q-m_j+1,\ldots,0\}$, then  $R_{m-q-1}(t)\equiv0$. This explains the condition stated in the lemma.

\smallskip

\noindent\textbf{Proof.} Claims (1) through (3) of the lemma will follow from the representation
$$
(\varepsilon-q)_q(\varepsilon+1)_{m-q-1}Q^{\varepsilon}_m(t)=(\varepsilon-q-t)_q(\varepsilon+1-t)_{m-q-1}+\frac{t}{\varepsilon}\left[A_1(t;\varepsilon)+{\varepsilon}A_2(t;\varepsilon)\right],
$$
where
\begin{multline*}
A_1(t;\varepsilon)=-(\varepsilon-q-t)_q\sum\limits_{k=0}^{m-q-1}(\bbbb)_{k}C_{k,r}(t)_{k}(\varepsilon+1-t)_{m-k-q-1}
\\
+\sum\limits_{k=m-q}^{m}(\bbbb)_{k}C_{k,r}(t+1)_{k-1}(\varepsilon-q-t)_{m-k}
\end{multline*}
and
$$
A_2(t;\varepsilon)=(\varepsilon-q-t)_q\sum\limits_{k=1}^{m-q-1}(\bbbb)_{k}C_{k,r}(t+1)_{k-1}(\varepsilon+1-t)_{m-k-q-1},
$$
which is verified by direct inspection.  Assume now that $\zeta(\varepsilon)$ is a root of $Q^{\varepsilon}_m(t)$, so that $Q^{\varepsilon}_m(\zeta(\varepsilon))=0$ or
$$
\zeta(\varepsilon)\left[A_1(\zeta(\varepsilon);\varepsilon)+{\varepsilon}A_2(\zeta(\varepsilon);\varepsilon)\right]=-\varepsilon(\varepsilon-q-\zeta(\varepsilon))_q(\varepsilon+1-\zeta(\varepsilon))_{m-q-1},
$$
since $(\varepsilon-q)_q(\varepsilon+1)_{m-q-1}$ does not vanish for sufficiently small $\varepsilon$.  The right hand side here goes to zero as  $\varepsilon\to0$
as long as $\zeta(\varepsilon)$ remains bounded. An alternative option $\zeta(\varepsilon)\to\infty$ must be discarded as the polynomial in $\zeta(\varepsilon)$ on the left hand side has degree $m$, while the polynomial on the right hand side has degree $m-1$, so that equality could not be maintained for sufficiently large $\zeta(\varepsilon)$.  Hence, the zeros of $Q^{\varepsilon}_m(t)$ tend to zeros of $tA_1(t;0)$. Clearly, $t=0$ is one of them. A more careful examination of
$$
A_1(t;0)=-(-q-t)_q\!\!\sum\limits_{k=0}^{m-q-1}(\bbbb)_{k}C_{k,r}(t)_{k}(1-t)_{m-k-q-1}
+\!\!\sum\limits_{k=m-q}^{m}(\bbbb)_{k}C_{k,r}(t+1)_{k-1}(-q-t)_{m-k}
$$
shows that
$$
A_1(t;0)=(-1)^{q+1}(t+1)_qR_{m-q-1}(t),
$$
where $R_{m-q-1}$ is defined in (\ref{eq:R}). This examination requires the relations $(-q-t)_q=(-1)^q(t+1)_q$,
\begin{multline*}
(t+1)_{k-1}(-q-t)_{m-k}=(t+1)_{q+k-m}(t+q+k-m+1)_{m-q-1}(-1)^{m-k}(t+q+k-m+1)_{m-k}
\\
=(-1)^{m-k}(t+1)_{q}(t+q+k-m+1)_{m-q-1}=(-1)^{k+q+1}(t+1)_{q}(1-t-k)_{m-q-1},
\end{multline*}
and
$$
(t)_{k}(1-t)_{m-k-q-1}=(-1)^{k}(1-t-k)_{m-q-1}.
$$
Hence, $t=-1,-2,\ldots,-q$ are the roots of $A_1(t;0)$. This establishes the claims (1) and (2).  Next, if $\zeta_1(\varepsilon)$ is the root tending zero as $\varepsilon\to0$, then
$$
\frac{\varepsilon}{\zeta_1(\varepsilon)}=-\frac{A_1(\zeta_1(\varepsilon);\varepsilon)+{\varepsilon}A_2(\zeta_1(\varepsilon);\varepsilon)}{(\varepsilon-q-\zeta_1(\varepsilon))_q(\varepsilon+1-\zeta_1(\varepsilon))_{m-q-1}}
\to-\frac{A_1(0;0)}{(-q)_q(1)_{m-q-1}}=\frac{R_{m-q-1}(0)}{(m-q-1)!},
$$
which proves the claim (3). Finally, to establish the claim (4), we will need the formula
$$
(a+b)_m=\sum_{k=0}^{m}\binom{m}{k}(a)_{k}(b)_{m-k}.
$$
Using this formula we get:
\begin{multline*}
(\varepsilon-q)_q(\varepsilon+1)_{m-q-1}Q^{\varepsilon}_m(-l)=\frac{1}{\varepsilon}\sum\limits_{k=0}^{m}(\bbbb)_{k}C_{k,r}(-l)_{k}(\varepsilon+l-q)_{m-k}
\\
=\frac{1}{\varepsilon}\sum\limits_{k=0}^{m}(\bbbb)_{k}C_{k,r}(-l)_{k}\sum\limits_{j=0}^{m-k}\binom{m-k}{j}(\varepsilon)_{j}(l-q)_{m-k-j}
\\
=\frac{1}{\varepsilon}\sum\limits_{k=0}^{l}(\bbbb)_{k}C_{k,r}(-l)_{k}(l-q)_{m-k}
+\frac{1}{\varepsilon}\sum\limits_{k=0}^{l}(\bbbb)_{k}C_{k,r}(-l)_{k}\sum\limits_{j=1}^{m-k}\binom{m-k}{j}(\varepsilon)_{j}(l-q)_{m-k-j}
\\
=\sum\limits_{k=0}^{l}(\bbbb)_{k}C_{k,r}(-l)_{k}\sum\limits_{j=1}^{m-k}\binom{m-k}{j}(\varepsilon+1)_{j-1}(l-q)_{m-k-j}
\\
\to\sum\limits_{k=0}^{l}(\bbbb)_{k}C_{k,r}(-l)_{k}\sum\limits_{j=1}^{m-k}\binom{m-k}{j}(j-1)!(l-q)_{m-k-j},
\end{multline*}
where we have taken account of
$$
(-l)_{k}(l-q)_{m-k}=0~\text{for}~l=0,\ldots,k~\text{and}~l=k+q-m+1,\ldots,q.
$$
As $q-m+1\le0$ this expression is guaranteed to vanish for $l=0,\ldots,q$.  Finally, using
$$
\sum\limits_{j=1}^{n}\binom{n}{j}(j-1)!(-r)_{n-j}=(-1)^rr!(n-r-1)!
$$
for $r=0,1,\ldots,n-1$, we have
$$
\sum\limits_{j=1}^{m-k}\binom{m-k}{j}(j-1)!(l-q)_{m-k-j}=(-1)^{q-l}(q-l)!(m-k-q+l-1)!
$$
which implies (\ref{eq:Qm0negative}). $\hfill\square$

\begin{lemma}\label{lm:hatQlimit}
Write $\hat{Q}_m^\varepsilon(t)$ for the polynomial \emph{(\ref{eq:Qmhat})} with $\cccc-\aaaa-m=-q+\varepsilon$, where $q\in\{0,1,\ldots,m-1\}$, $\varepsilon>0$
and let $\eta_1(\varepsilon),\ldots,\eta_m(\varepsilon)$ denote its zeros.   Suppose $f_j-\aaaa\notin\{m-q-m_j,m-q-m_j+1,\ldots,0\}$ for indices $0\le{j}\le{r}$ such that $m-q-m_j\le0$. Then the next claims are true as $\varepsilon\to0$\emph{:}

\smallskip

\emph{(1)} the first \emph{(}after possible renumbering\emph{)} $q+1$ zeros tend to consecutive non-positive integers\emph{:} $\eta_1\to0$, $\eta_2\to-1$, $\ldots$, $\eta_{q+1}\to-q$\emph{;}

\smallskip

\emph{(2)} the remaining $m-q-1$ zeros \emph{(}no such zeros remain if $q=m-1$\emph{)} tend to the zeros of the reduced polynomial
\begin{multline}\label{eq:hatR}
\hat{R}_{m-q-1}(t)=\frac{(\bbbb)_{q+1}}{(\bbbb-\aaaa)_{q+1}}\sum\limits_{k=0}^{q}\frac{(-m+k)_{q-k+1}(\aaaa)_{k}C_{k,r}}{(q-k+1)!}
{}_{3}F_{2}\!\!\left(\begin{matrix}1-m+q,t+q+1,1-\bbbb-k\\\aaaa-\bbbb+1,2-k+q\end{matrix}\right)
\\
+\sum\limits_{k=q+1}^{m}\frac{(-1)^k(\aaaa)_{k}(\bbbb)_{k}C_{k,r}(t+q+1)_{k-q-1}}{(\aaaa-\bbbb-q)_k(k-q-1)!}
{}_{3}F_{2}\!\!\left(\begin{matrix}-m+k,t+k,-\bbbb-q\\\aaaa-\bbbb-q+k,k-q\end{matrix}\right)
\end{multline}
of degree $m-q-1$\emph{;}

\smallskip

\emph{(3)} the limit relation
\begin{equation}\label{eq:hatA00}
\lim\limits_{\varepsilon\to0}\frac{\varepsilon}{\eta_1(\varepsilon)}=(-1)^{q+1}\hat{R}_{m-q-1}(0)
\end{equation}
holds\emph{;}

\smallskip

\emph{(4)}  the values $\hat{Q}_m^\varepsilon(-l)$ at first $q+1$ non-positive integers converge  to the following expressions\emph{:}
\begin{equation}\label{eq:hatQm0negative}
\hat{Q}_m^0(-l)=\sum\limits_{k=0}^{l}\frac{(-1)^k(-l)_{k}(\aaaa)_k(\bbbb)_kC_{k,r}}{(-q)_k(\aaaa-\bbbb-q)_k}
{}_{3}F_{2}\!\!\left(\begin{matrix}k-m,k-l,-\bbbb-q\\k-q,-\bbbb-q+\aaaa+k\end{matrix}\right),~~l=1,\ldots,q,
\end{equation}
and $\hat{Q}_m^0(0)=1$.
\end{lemma}
\textbf{Remark.}  We claim that if $m-q-m_j\le0$ for some $0\le{j}\le{r}$ and $f_j-b\in\{m-q-m_j,m-q-m_j+1,\ldots,0\}$, then $\hat{R}_{m-q-1}(t)\equiv0$ just like
$R_{m-q-1}(t)\equiv0$. Unfortunately we do not have an explicit formula that would make this claim obvious. However, we can argue as follows.  The transformation (\ref{eq:KRPTh2limit}) can be obtained by an application of (\ref{eq:KRPTh1-1}) to the hypergeometric function on the right hand side of (\ref{eq:KRPTh1limit}).
This means that the polynomial $\hat{R}_{m-q-1}(t)$ is proportional to the polynomial $Q_m(\bbbb,\cccc,\llambda+q+1|t)$, where $Q_m$ is defined in (\ref{eq:Qm})
with $\m=(1,\ldots,1)$.  But for this $\m$ the coefficients of $Q_m$ are symmetric functions of $\llambda+q+1$ and so of $\llambda$. As $\llambda$ are the roots of
$R_{m-q-1}(t)$ we conclude that the coefficients of $Q_m(\bbbb,\cccc,\llambda+q+1|t)$ and hence of $\hat{R}_{m-q-1}(t)$  can be written in terms of coefficients
of $R_{m-q-1}(t)$.  This explains why these coefficients vanish if the coefficients of $\hat{R}_{m-q-1}(t)$ vanish.  The details will be given elsewhere.

\textbf{Proof.} All claims of the lemma will follow from the representation
\begin{equation}\label{eq:hatQmdecomp}
\hat{Q}_m^\varepsilon(t)=1+\frac{t}{\varepsilon}\left[A_1(t;\varepsilon)+{\varepsilon}A_2(t;\varepsilon)\right]
\end{equation}
with
\begin{multline*}
A_1(t;\varepsilon)=\frac{(t+1)_q}{(\varepsilon-q)_{q}}\biggl\{
\sum\limits_{k=0}^{q}\frac{(-1)^kC_{k,r}(\aaaa)_k(\bbbb)_k(-m+k)_{q-k+1}(\varepsilon-\bbbb-q)_{q-k+1}}{(\aaaa-\bbbb-q+\varepsilon)_{q+1}(q-k+1)!}
\\
\times{}_{4}F_{3}\!\!\left(\begin{matrix}-m+q+1,t+q+1,\varepsilon-\bbbb-k+1,1\\\varepsilon+1, \aaaa-\bbbb+1+\varepsilon,2+q-k\end{matrix}\right)
+\sum\limits_{k=q+1}^{m}\frac{(-1)^kC_{k,r}(\aaaa)_k(\bbbb)_k(t+q+1)_{k-q-1}}{(\varepsilon+1)_{k-q-1}(\aaaa-\bbbb-q+\varepsilon)_{k}}
\\
\times{}_{3}F_{2}\!\!\left(\begin{matrix}-m+k,t+k,\varepsilon-\bbbb-q\\\aaaa-\bbbb-q+k+\varepsilon,\varepsilon-q+k\end{matrix}\right)
\biggr\}
\end{multline*}
and
\begin{multline*}
A_2(t;\varepsilon)=\sum\limits_{j=1}^{q}\frac{(-m)_j(\varepsilon-\bbbb-q)_j}{(\varepsilon-q)_j(\aaaa-\bbbb-q+\varepsilon)_jj!}
\\
+\sum\limits_{k=1}^{q}\frac{(-1)^kC_{k,r}(\aaaa)_k(\bbbb)_k(t+1)_{k-1}}{(\varepsilon-q)_k(\aaaa-\bbbb-q+\varepsilon)_k}
\sum\limits_{j=0}^{q-k}
\frac{(-m+k)_j(t+k)_j(\varepsilon-\bbbb-q)_j}{(\varepsilon-q+k)_j(\aaaa-\bbbb-q+k+\varepsilon)_jj!}.
\end{multline*}
The proof of  decomposition (\ref{eq:hatQmdecomp}) is by substituting $A_1$, $A_2$ into (\ref{eq:hatQmdecomp}) and tedious but elementary calculation.
Note that $A_1(0;0)$ does not vanish.  Hence, the equation $\hat{Q}_m^\varepsilon(t)=0$ or
$$
t\left[A_1(t;\varepsilon)+{\varepsilon}A_2(t;\varepsilon)\right]=-\varepsilon
$$
implies that $t\to0$ as $\varepsilon\to0$, so that one root, call it $\eta_1(\varepsilon)$ tends to zero. Further, the second term in brackets vanishes so that the remaining roots
converge to the roots of $A_1(t;0)$.  The factor $(t+1)_q$ shows that $-1,-2,\ldots,-q$ are zeros of $A_1(t;0)$ which proves the claim (1).  Further, noting that
$$
A_1(t;0)=\frac{(t+1)_{q}}{(-q)_q}\hat{R}_{m-q-1}(t),
$$
where $\hat{R}_{m-q-1}(t)$ is defined in (\ref{eq:hatR}), we establish the second claim.  Relation (\ref{eq:hatA00}) again follows from (\ref{eq:hatQmdecomp}). Indeed,
$$
\frac{\varepsilon}{\eta_1(\varepsilon)}=-A_1(t;\varepsilon)-{\varepsilon}A_2(t;\varepsilon)\to -A_1(0;0)=-\frac{(1)_{q}}{(-q)_q}\hat{R}_{m-q-1}(0)=(-1)^{q+1}\hat{R}_{m-q-1}(0)
$$
for $\eta_1(\varepsilon)\to0$ as $\varepsilon\to0$.  Finally, as $A_1(-l;\varepsilon)=0$ for $l=1,\ldots,q$ we get from (\ref{eq:hatQmdecomp})
$$
\hat{Q}_m^\varepsilon(-l)=1-lA_2(-l;\varepsilon)\to1-lA_2(-l;0)~\text{as}~\varepsilon\to0.
$$
This can be seen to equal (\ref{eq:hatQm0negative}) by manipulating $A_2(-l;0)$. Note that the straightforward substitution $t=-l$ into (\ref{eq:Qmhat}) and letting $\varepsilon\to0$ gives similar expression but with summation running up to $k=m$ and thus containing many undefined $0/0$ terms. $\hfill\square$

Our main results are presented in the following three theorems.

\begin{theorem}\label{th:KRPTh1limit}
Suppose $q\in\{0,\ldots,m-1\}$,   $\bbbb+m-q\notin-\N_0$ and  
$f_j-\bbbb\notin\{m-q-m_j,m-q-m_j+1,\ldots,0\}$ for indices $0\le{j}\le{r}$ such that $m-q-m_j\le0$. Then for all $x\in\C\!\setminus\![1,\infty)$ the following identity holds\emph{:}
\begin{multline}\label{eq:KRPTh1limit}
(1-x)^{\aaaa}{}_{r+2}F_{r+1}\left.\!\left(\!\begin{matrix}\aaaa,\bbbb,\f+\m\\\bbbb+m-q, \f\end{matrix}\right\vert x\!\right)=
\sum\limits_{j=0}^{q}\frac{(\aaaa)_j(-q)_jQ^{0}_m(-j)}{(\bbbb+m-q)_jj!}\left(\frac{x}{x-1}\right)^j
\\
+\frac{x^{q+1}(\aaaa)_{q+1}R_{m-q-1}(-q-1)}{(x-1)^{q+1}(\bbbb+m-q)_{q+1}(m-q-1)!}
{}_{m-q+1}F_{m-q}\left.\!\!\left(\begin{matrix}1,\aaaa+q+1,\llambda+q+2\\\bbbb+m+1, \llambda+q+1\end{matrix}\right\vert\frac{x}{x-1}\!\right),
\end{multline}
where $\llambda=(\lambda_1,\ldots,\lambda_{m-q-1})$ are the zeros of the polynomial $R_{m-q-1}(t)$ of degree $m-q-1$ defined in \emph{(\ref{eq:R})}, and
the values $Q^0_m(-j)$  \emph{(}$0\le{j}\le{q}$\emph{)} are given in \emph{(\ref{eq:Qm0negative})}.  The number $R_{m-q-1}(-q-1)$ is computed by
\begin{equation}\label{eq:R-q-1}
R_{m-q-1}(-q-1)=(-1)^{m-q+1}\sum\limits_{k=0}^{q+1}(\bbbb)_{k}C_{k,r}(-1)^{k}(k-m)_{m-q-1},
\end{equation}
where $C_{k,r}$ is given by \emph{(\ref{eq:Ckr})}.
\end{theorem}
\textbf{Proof.} Writing $\cccc-\bbbb=m-q+\varepsilon$ and $y=x/(x-1)$ formula (\ref{eq:KRPTh1-1}) takes the form
\begin{multline}\label{eq:Th1limit-proof}
(1-x)^{\aaaa}{}_{r+2}F_{r+1}\left.\!\!\!\left(\begin{matrix}\aaaa, \bbbb, \f+\m\\\bbbb+m-q+\varepsilon,\f\end{matrix}\right\vert x\right)
=\sum\limits_{j=0}^{q}\frac{(\aaaa)_j(\varepsilon-q)_j(\zetta(\varepsilon)+1)_j}{(\bbbb+m-q+\varepsilon)_j(\zetta(\varepsilon))_jj!}y^j
\\
+\sum\limits_{j=q+1}^{\infty}\frac{(\aaaa)_j(\varepsilon-q)_{j}(\zetta(\varepsilon)+1)_{j}}{(\bbbb+m-q+\varepsilon)_{j}(\zetta(\varepsilon))_{j}j!}y^j
=1+\sum\limits_{j=1}^{q}\frac{(\aaaa)_j(\varepsilon-q)_j}{(\bbbb+m-q+\varepsilon)_jj!}\frac{Q^{\varepsilon}_m(-j)}{Q^{\varepsilon}_m(0)}y^j
\\
+\frac{\varepsilon}{\zeta_1(\varepsilon)}\sum\limits_{j=q+1}^{\infty}\frac{(\aaaa)_j(\varepsilon-q)_q(\varepsilon+1)_{j-q-1}(\zetta(\varepsilon)+j)}{(\bbbb+m-q+\varepsilon)_j(\zetta_{[1]}(\varepsilon))j!}y^{j},
\end{multline}
where we applied the straightforward relation
$$
\frac{(\zetta(\varepsilon)+1)_j}{(\zetta(\varepsilon))_j}=\frac{(\zetta(\varepsilon)+j)_1}{(\zetta(\varepsilon))_1}=\frac{Q^{\varepsilon}_m(-j)}{Q^{\varepsilon}_m(0)}
$$
(recall that $Q_m^\varepsilon(t)$ is the polynomial (\ref{eq:Qm}) with $\cccc-\bbbb-m=-q+\varepsilon$).   Next, denoting by $\llambda=(\lambda_1,\ldots,\lambda_{m-q-1})$ the zeros of the reduced polynomial $R_{m-q-1}$ defined in (\ref{eq:R}) we compute using claims (1) and (2) of Lemma~\ref{lm:Qlimit}:
\begin{multline*}
\lim\limits_{\varepsilon\to0}\frac{(\varepsilon-q)_q(\varepsilon+1)_{j-q-1}(\zetta(\varepsilon)+j)}{(\bbbb+m-q+\varepsilon)_j(\zetta_{[1]}(\varepsilon))j!}=
\frac{(-q)_q(1)_{j-q-1}j(j-1)\cdots(j-q)(\llambda+j)_1}{(\bbbb+m-q)_j(-1)(-2)\cdots(-q)(\llambda)_1j!}
\\
=\frac{(\llambda+j)_1}{(\bbbb+m-q)_j(\llambda)_1}=\frac{(\llambda+q+2)_{j-q-1}R_{m-q-1}(-q-1)}{(\bbbb+m-q)_{q+1}(\bbbb+m+1)_{j-q-1}(\llambda+q+1)_{j-q-1}R_{m-q-1}(0)}.
\end{multline*}
Substituting this limit in (\ref{eq:Th1limit-proof}) and invoking (\ref{eq:A00Q}) we arrive at
\begin{multline*}
(1-x)^{\aaaa}{}_{r+2}F_{r+1}\left.\!\!\!\left(\begin{matrix}\aaaa, \bbbb, \f+\m\\\bbbb+m-q,\f\end{matrix}\right\vert x\right)
=1+\sum\limits_{j=1}^{q}\frac{(\aaaa)_j(-q)_j}{(\bbbb+m-q)_jj!}\frac{Q^{0}_m(-j)}{Q^{0}_m(0)}y^j
\\
+\frac{(\aaaa)_{q+1}R_{m-q-1}(0)R_{m-q-1}(-q-1)}{R_{m-q-1}(0)(\bbbb+m-q)_{q+1}(m-q-1)!}
\sum\limits_{j=q+1}^{\infty}\frac{(\aaaa+q+1)_{j-q-1}(\llambda+q+2)_{j-q-1}}{(\bbbb+m+1)_{j-q-1}(\llambda+q+1)_{j-q-1}}y^{j}.
\end{multline*}
In view of $Q^{0}_m(0)=1$ by Lemma~\ref{lm:Qlimit}(4), this coincides with (\ref{eq:KRPTh1limit}). Finally, formula (\ref{eq:R-q-1}) follows directly from (\ref{eq:R}). $\hfill\square$

\medskip

The Kummer-type transformation (\ref{eq:Kummer-type}) degenerates as follows.

\begin{theorem}\label{th:Kummerlimit}
Suppose $q\in\{0,\ldots,m-1\}$,   $\bbbb+m-q\notin-\N_0$ and
$f_j-\bbbb\notin\{m-q-m_j,m-q-m_j+1,\ldots,0\}$ for indices $0\le{j}\le{r}$ such that $m-q-m_j\le0$. Then for all $x\in\C$
\begin{multline}\label{eq:Kummerlimit}
e^{-x}{}_{r+1}F_{r+1}\left.\!\left(\!\begin{matrix}\bbbb,\f+\m\\\bbbb+m-q, \f\end{matrix}\right\vert x\!\right)=\sum\limits_{j=0}^{q}\frac{(-q)_jQ^{0}_m(-j)}{(\bbbb+m-q)_jj!}(-x)^j
\\
+\frac{(-x)^{q+1}R_{m-q-1}(-q-1)}{(\bbbb+m-q)_{q+1}(m-q-1)!}
{}_{m-q}F_{m-q}\left.\!\!\left(\begin{matrix}1,\llambda+q+2\\\bbbb+m+1, \llambda+q+1\end{matrix}\right\vert-x\!\right),
\end{multline}
where  $\llambda=(\lambda_1,\ldots,\lambda_{m-q-1})$ are zeros of the polynomial $R_{m-q-1}(t)$ defined in \emph{(\ref{eq:R})},
the values $Q^0_m(-j)$ \emph{(}$0\le{j}\le{q}$\emph{)} are given in \emph{(\ref{eq:Qm0negative})}, and the number $R_{m-q-1}(-q-1)$ is found in \emph{(\ref{eq:R-q-1})}.
\end{theorem}

The proof of this theorem repeats verbatim the proof of Theorem~\ref{th:KRPTh1limit} and thus will be omitted.

\begin{theorem}\label{th:KRPTh2limit}
Suppose $q\in\{0,\ldots,m-1\}$, $\aaaa+m-q\notin-\N_0$, $\aaaa-\bbbb\notin\{q+1-m,\ldots,q\}$ and $f_j-\aaaa\notin\{m-q-m_j,m-q-m_j+1,\ldots,0\}$ for indices $0\le{j}\le{r}$ such that $m-q-m_j\le0$. Then for all $x\in\C\!\setminus\![1,\infty)$ the following identity holds\emph{:}
\begin{multline}\label{eq:KRPTh2limit}
(1-x)^{\bbbb+q}{}_{r+2}F_{r+1}\left.\!\!\left(\begin{matrix}\aaaa,\bbbb,\f+\m\\\aaaa+m-q, \f\end{matrix}\right\vert x\!\right)=\sum\limits_{j=0}^{q}\frac{(-q)_j(\aaaa-\bbbb-q)_j}{(\aaaa+m-q)_jj!}\hat{Q}^0_m(-j)x^j
\\
+x^{q+1}\frac{\hat{R}_{m-q-1}(-q-1)(\bbbb-\aaaa)_{q+1}}{(\aaaa+m-q)_{q+1}}
{}_{m-q+1}F_{m-q}\left.\!\!\left(\begin{matrix}1,\aaaa-\bbbb+1,\ggamma+q+2\\\aaaa+m+1, \ggamma+q+1\end{matrix}\right\vert x\!\right),
\end{multline}
where $\ggamma=(\gamma_1,\ldots,\gamma_{m-q-1})$ are the zeros of the polynomial $\hat{R}_{m-q-1}(t)$ of degree $m-q-1$ defined in \emph{(\ref{eq:hatR})} and
the values $\hat{Q}^0_m(-j)$, $j=0,\ldots,q$ are given in \emph{(\ref{eq:hatQm0negative})}.  The number $\hat{R}_{m-q-1}(-q-1)$ is computed by
\begin{equation}\label{eq:hatR-q-1}
\hat{R}_{m-q-1}(-q-1)=\frac{(\bbbb)_{q+1}}{(\bbbb-\aaaa)_{q+1}}\sum\limits_{k=0}^{q+1}\frac{(-m+k)_{q+1-k}(\aaaa)_kC_{k,r}}{(q+1-k)!},
\end{equation}
where $C_{k,r}$ is given by \emph{(\ref{eq:Ckr})}.
\end{theorem}
\textbf{Proof.} Writing $\cccc-\aaaa=m-q+\varepsilon$ formula (\ref{eq:KRPTh1-2}) takes the form
\begin{multline}\label{eq:Th2limit-proof}
(1-x)^{\bbbb+q-\varepsilon}{}_{r+2}F_{r+1}\left.\!\!\!\left(\begin{matrix}\aaaa, \bbbb, \f+\m\\\aaaa+m-q+\varepsilon,\f\end{matrix}\right\vert x\right)
=\sum\limits_{j=0}^{q}\frac{(\varepsilon-q)_j(\aaaa-\bbbb-q+\varepsilon)_j (\etta(\varepsilon)+1)_j}{(\aaaa+m-q+\varepsilon)_j(\etta(\varepsilon))_jj!}x^j
\\
+\sum\limits_{j=q+1}^{\infty}\frac{(\varepsilon-q)_q\varepsilon(\varepsilon+1)_{j-q-1}(\aaaa-\bbbb-q+\varepsilon)_j(\etta(\varepsilon)+1)_j}{(\aaaa+m-q+\varepsilon)_j(\etta(\varepsilon))_jj!}x^j
\\
=1+\sum\limits_{j=1}^{q}\frac{(\varepsilon-q)_j(\aaaa-\bbbb-q+\varepsilon)_j}{(\aaaa+m-q+\varepsilon)_jj!}\frac{\hat{Q}^{\varepsilon}_m(-j)}{\hat{Q}^{\varepsilon}_m(0)}x^j
\\
+\frac{\varepsilon}{\eta_1(\varepsilon)}\sum\limits_{j=q+1}^{\infty}\frac{(\varepsilon-q)_q(\varepsilon+1)_{j-q-1}(\aaaa-\bbbb-q+\varepsilon)_j(\etta(\varepsilon)+j)_1}{(\aaaa+m-q+\varepsilon)_j(\etta_{[1]}(\varepsilon))_1j!}x^j,
\end{multline}
where we used the straightforward fact that
$$
\frac{(\etta(\varepsilon)+1)_j}{(\etta(\varepsilon))_j}=\frac{(\etta(\varepsilon)+j)_1}{(\etta(\varepsilon))_1}=\frac{\hat{Q}^{\varepsilon}_m(-j)}{\hat{Q}^{\varepsilon}_m(0)}
$$
(recall that $\hat{Q}_m^\varepsilon(t)$ is the polynomial (\ref{eq:Qmhat}) with $\cccc-\aaaa-m=-q+\varepsilon$). Letting $\varepsilon\to0$ we see that the first sum tends to the first term on the right hand side of (\ref{eq:KRPTh2limit}) in view of  $\hat{Q}^{\varepsilon}_m(0)=1$ by Lemma~\ref{lm:hatQlimit}(4). To find the limit of the second sum we resort to Lemma~\ref{lm:hatQlimit}. First, the factor $\varepsilon/\eta_1(\varepsilon)$ tends to $(-1)^{q+1}\hat{R}_{m-q-1}(0)$ by (\ref{eq:hatA00}). Next, by claims (1) and (2) of Lemma~\ref{lm:hatQlimit} we have
$$
(\etta(\varepsilon)+j)_1\to j(j-1)\cdots(j-q)(\lambda_1+j)\cdots(\lambda_{m-q-1}+j),
$$
where $\lambda_1,\ldots,\lambda_{m-q-1}$ are the roots of $\hat{R}_{m-q-1}(t)$ defined in (\ref{eq:hatR}). Similar formula holds for $j=0$.  Then for $j\ge{q+1}$ as $\varepsilon\to0$ \begin{multline*}
\frac{(\varepsilon-q)_q(\varepsilon+1)_{j-q-1}(\aaaa-\bbbb-q+\varepsilon)_{j}(\etta(\varepsilon)+j)_1}{(\aaaa+m-q+\varepsilon)_{j}(\etta_{[1]}(\varepsilon))_1j!}x^j
\\
\to\frac{(-q)_q(1)_{j-q-1}(\aaaa-\bbbb-q)_{j}j(j-1)\cdots(j-q)(\ggamma+j)_1}{(\aaaa+m-q)_j (-1)\cdots(-q)(\ggamma)_1j!}x^j
\\
=\frac{(\aaaa-\bbbb-q)_{q+1}(\aaaa-\bbbb+1)_{j-q-1}(\ggamma+1)_{q+1}(\ggamma+q+2)_{j-q-1}}{(\aaaa+m-q)_{q+1}(\aaaa+m+1)_{j-q-1}(\ggamma)_{q+1}(\ggamma+q+1)_{j-q-1}}x^j
\\
=\frac{\hat{R}_{m-q-1}(-q-1)}{\hat{R}_{m-q-1}(0)}\frac{(-1)^{q+1}(\bbbb-\aaaa)_{q+1}(\aaaa-\bbbb+1)_{j-q-1}(\ggamma+q+2)_{j-q-1}}{(\aaaa+m-q)_{q+1}(\aaaa+m+1)_{j-q-1}(\ggamma+q+1)_{j-q-1}}x^j.
\end{multline*}
So letting $\varepsilon\to0$ in  (\ref{eq:Th2limit-proof}) and substituting these limits we arrive at (\ref{eq:KRPTh2limit}). It remains to compute $\hat{R}_{m-q-1}(-q-1)$ in the form (\ref{eq:hatR-q-1}). To this end just substitute $t=-q-1$ into (\ref{eq:hatR}) and notice that the second sum reduces to $q+1$-th term of the first sum. $\hfill\square$

\section{Extreme cases}

In this section we consider the cases $q=0$, $q=m-1$ and $q=m-2$ of the general theorems from the previous section.
The case $q=0$ leads to noticeable simplification of the terms outside the generalized hypergeometric functions on the right hand sides.
For $q=m-1$ the characteristic polynomials $R_{m-q-1}$ and $\hat{R}_{m-q-1}$ degenerate to constants, while the generalized hypergeometric functions
on the right hand sides reduce to the Gauss hypergeometric functions ${}_2F_{1}$.  Finally, when $q=m-2$ the characteristic polynomials $R_{m-q-1}$ and $\hat{R}_{m-q-1}$
become linear and the root can be explicitly computed, while  the generalized hypergeometric functions on the right hand sides are Clausen's functions ${}_3F_{2}$.
The cases $q=0$ of Theorem~\ref{th:KRPTh1limit} and Theorem~\ref{th:Kummerlimit} take the following form (recall that $(\f)_{\m}=(f_1)_{m_1}\cdots(f_r)_{m_r}$).
\begin{corollary}\label{cr:KRPTh1limit0}
Suppose $\bbbb+m\notin-\N_0$ and $f\ne{\bbbb}$ if $m=r=1$, then the following identities hold\emph{:}
\begin{equation}\label{eq:KRPTh1-1-reduced}
(1-x)^{\aaaa}{}_{r+2}F_{r+1}\left.\!\left(\!\begin{matrix}\aaaa, \bbbb, \f+\m\\\bbbb+m,\f\end{matrix}\right\vert x\!\right)
=\frac{(\bbbb)_m}{(\f)_{\m}}+\biggl(1-\frac{(\bbbb)_m}{(\f)_{\m}}\biggr){}_{m+1}F_{m}\left.\!\left(\begin{matrix}1, \aaaa, \llambda+1\\\bbbb+m, \llambda\end{matrix}\right\vert\frac{x}{x-1}\!\right)
\end{equation}
for all $x\in\C\!\setminus\![1,\infty)$  and
\begin{equation}\label{eq:KRPThKummerq-0}
e^{-x}{}_{r+1}F_{r+1}\left.\!\left(\begin{matrix}\bbbb,\f+\m\\\bbbb+m, \f\end{matrix}\right\vert x\!\right)=
\frac{(\bbbb)_m}{(\f)_{\m}}+\biggl(1-\frac{(\bbbb)_m}{(\f)_{\m}}\biggr){}_{m}F_{m}\left.\!\!\left(\begin{matrix}1,\llambda+1\\\bbbb+m, \llambda\end{matrix}\right\vert-x\!\right)
\end{equation}
for all $x\in\C$, where $\llambda=(\lambda_1,\ldots,\lambda_{m-1})$ are the roots of the polynomial
\begin{equation}\label{eq:Rm-1}
R_{m-1}(t)=\sum\limits_{k=0}^{m}(\bbbb)_k(-1)^kC_{k,r}(1-t-k)_{m-1}.
\end{equation}
\end{corollary}
\textbf{Proof.}  Indeed, taking  $q=0$ in (\ref{eq:KRPTh1limit}) we obtain:
$$
(1-x)^{\aaaa}{}_{r+2}F_{r+1}\left.\!\left(\!\begin{matrix}\aaaa,\bbbb,\f+\m\\\bbbb+m, \f\end{matrix}\right\vert x\!\right)=
1+\frac{x{\aaaa}R_{m-1}(-1)}{(x-1)(\bbbb+m)(m-1)!}
{}_{m+1}F_{m}\left.\!\!\left(\begin{matrix}1,\aaaa+1,\llambda+2\\\bbbb+m+1, \llambda+1\end{matrix}\right\vert\frac{x}{x-1}\!\right),
$$
where we used that $Q_m^0(0)=1$.  Further, since
$$
R_{m-1}(-1)\!=\!\frac{R_{m-1}(-1)}{R_{m-1}(0)}R_{m-1}(0)\!=\!\frac{(\llambda+1)_1}{(\llambda)_1}(m-1)!\left(1-C_{m,r}(\bbbb)_m\right)\!
=\!\!\biggl[1-\frac{(\bbbb)_m}{(\f)_{\m}}\biggr]\!\frac{(\llambda+1)_1(m-1)!}{(\llambda)_1}
$$
by definition of $R_{m-1}(t)$ and the expression for $C_{m,r}$ below (\ref{eq:Qm}).  Then,
\begin{multline*}
\frac{x{\aaaa}R_{m-1}(-1)}{(x-1)(\bbbb+m)(m-1)!}
{}_{m+1}F_{m}\left.\!\!\left(\begin{matrix}1,\aaaa+1,\llambda+2\\\bbbb+m+1, \llambda+1\end{matrix}\right\vert\frac{x}{x-1}\!\right)
\\
=\left(1-\frac{(\bbbb)_m}{(\f)_{\m}}\right)\frac{x{\aaaa}(\llambda+1)_1}{(x-1)(\bbbb+m)(\llambda)_1}
{}_{m+1}F_{m}\left.\!\!\left(\begin{matrix}1,\aaaa+1,\llambda+2\\\bbbb+m+1, \llambda+1\end{matrix}\right\vert\frac{x}{x-1}\!\right)
\\
=\left(1-\frac{(\bbbb)_m}{(\f)_{\m}}\right)\biggl[
{}_{m+1}F_{m}\left.\!\!\left(\begin{matrix}1,\aaaa,\llambda+1\\\bbbb+m, \llambda\end{matrix}\right\vert\frac{x}{x-1}\!\right)-1\biggr],
\end{multline*}
which yields (\ref{eq:KRPTh1-1-reduced}).
Identity (\ref{eq:KRPThKummerq-0}) can be proved starting from (\ref{eq:Kummerlimit}) in a similar fashion or by applying
confluence to (\ref{eq:KRPTh1-1-reduced}). $\hfill\square$

Next we present the $q=0$ case  of Theorem~\ref{th:KRPTh2limit}.
\begin{corollary}\label{cr:KRPTh2limit0}
If $\aaaa+m\notin-\N_0$, $\aaaa-\bbbb\notin\{1-m,\ldots,0\}$ and $f\ne{\aaaa}$ if $m=r=1$, then for all $x\in\C\!\setminus\![1,\infty)$ the following identity holds\emph{:}
\begin{equation}\label{eq:KRPTh1-2-reduced}
(1-x)^{\bbbb}{}_{r+2}F_{r+1}\left.\!\left(\!\begin{matrix}\aaaa,\bbbb,\f+\m\\\aaaa+m, \f\end{matrix}\right\vert x\!\right)
=1
+x\bbbb\biggl(\frac{\aaaa(\f+\m)}{(\aaaa+m)(\f)}-1\biggr)
{}_{m+1}F_{m}\left.\!\left(\!\begin{matrix}1,\aaaa-\bbbb+1,\ggamma+2\\\aaaa+m+1, \ggamma+1\end{matrix}\right\vert x\!\right),
\end{equation}
\end{corollary}
where $\ggamma=(\gamma_1,\ldots,\gamma_{m-1})$ are the roots of the polynomial
\begin{multline}\label{eq:hatRm-1}
\hat{R}_{m-1}(t)=\frac{(-m)(-\bbbb)C_{0,r}}{(\aaaa-\bbbb)}{}_{3}F_{2}\!\left(\begin{matrix}-m+1,t+1,1-\bbbb\\2,\aaaa-\bbbb+1\end{matrix}\right)
\\
+\sum\limits_{k=1}^{m}\frac{(-1)^kC_{k,r}(\aaaa)_k(\bbbb)_k(t+1)_{k-1}}{(1)_{k-1}(\aaaa-\bbbb)_k}{}_{3}F_{2}\!\left(\begin{matrix}-m+k,t+k,-\bbbb\\k,\aaaa-\bbbb+k\end{matrix}\right).
\end{multline}
\textbf{Proof.}  Indeed, taking  $q=0$ in (\ref{eq:KRPTh2limit}) and using the value of $\hat{R}_{m-1}(-1)$ from (\ref{eq:hatR-q-1})
we get (\ref{eq:KRPTh1-2-reduced}) after rather straightforward simplifications.  $\hfill\square$

\textbf{Remark.}  If we apply the procedure similar to the one used in the proof of Corollary~\ref{cr:KRPTh1limit0} we get a somewhat more complicated formula:
\begin{equation}\label{eq:KRPTh1-2-reducedAlt}
(1-x)^{\bbbb}{}_{r+2}F_{r+1}\left.\!\left(\!\begin{matrix}\aaaa,\bbbb,\f+\m\\\aaaa+m, \f\end{matrix}\right\vert x\!\right)
=A+(1-A){}_{m+1}F_{m}\left.\!\left(\!\begin{matrix}1, \aaaa-\bbbb, \ggamma+1\\\aaaa+m, \ggamma\end{matrix}\right\vert x\!\right),
\end{equation}
where $\aaaa-\bbbb\notin\{-m+1,-m+2,\ldots,0\}$,
$$
A=\sum\limits_{k=0}^{m}\frac{(-1)^kC_{k,r}(\aaaa)_k(\bbbb)_k(1-\aaaa-m)_{m-k}}{(\aaaa-\bbbb)_k(1-\cccc+\bbbb)_{m-k}}.
$$

Another extreme case $q=m-1$ is treated in the next two corollaries.
\begin{corollary}\label{cr:KRPTh1limitm-1}
If $\bbbb+1\notin-\N_0$ and $f_j-\bbbb\notin\{1-m_j,2-m_j,\ldots,0\}$ \emph{(}$1\le{j}\le{r}$\emph{)}, then the following identities hold\emph{:}
\begin{multline}\label{eq:KRPTh1limitm-1}
(1-x)^{\aaaa}{}_{r+2}F_{r+1}\left.\!\left(\begin{matrix}\aaaa,\bbbb,\f+\m\\\bbbb+1, \f\end{matrix}\right\vert x\!\right)=
\sum\limits_{j=0}^{m-1}\frac{(\aaaa)_j}{(\bbbb+1)_{j}}\left(\frac{x}{x-1}\right)^j\sum\limits_{k=0}^{j}(-1)^k(\bbbb)_{k}C_{k,r}
\\
+\frac{x^{m}(\aaaa)_{m}(\f-\bbbb)_{\m}}{(x-1)^{m}(\bbbb+1)_{m}(\f)_{\m}}
{}_{2}F_{1}\left.\!\!\left(\begin{matrix}1, \aaaa+m\\\bbbb+m+1 \end{matrix}\right\vert\frac{x}{x-1}\!\right)
\end{multline}
for $x\in\C\!\setminus\![1,\infty)$ and
\begin{multline}\label{eq:KRPThKummerq-m-1}
e^{-x}{}_{r+1}F_{r+1}\left.\!\left(\begin{matrix}\bbbb,\f+\m\\\bbbb+1, \f\end{matrix}\right\vert x\!\right)=
\sum\limits_{j=0}^{m-1}\frac{(-x)^j}{(\bbbb+1)_{j}}\sum\limits_{k=0}^{j}(-1)^k(\bbbb)_{k}C_{k,r}
\\
+\frac{(-x)^{m}(\f-\bbbb)_{\m}}{(\bbbb+1)_{m}(\f)_{\m}}
{}_{1}F_{1}\left.\!\!\left(\begin{matrix}1\\\bbbb+m+1 \end{matrix}\right\vert-x\!\right)
\end{multline}
for all $x\in\C$.  The numbers $C_{k,r}$ are defined in \emph{(\ref{eq:Ckr})}.
\end{corollary}
\textbf{Proof.}  For $q=m-1$ the polynomial $R_{m-q-1}=R_0$ defined in (\ref{eq:R}) reduces to the  constant
$$
R_{0}=\sum\limits_{k=0}^{m}(\bbbb)_{k}(-1)^{k}C_{k,r}=\sum\limits_{k=0}^{m}\frac{(\bbbb)_{k}}{k!}{}_{r+1}F_{r}\!\left(\begin{matrix}-k,\f+\m\\\f\end{matrix}\right)
=\frac{(\f-\bbbb)_{\m}}{(\f)_{\m}},
$$
where the last equality is \cite[Corollary~3]{MP2012MathComm}.  Further, from (\ref{eq:Qm0negative})
$$
Q_m^0(-j)=\frac{1}{(1-m)_j}\sum\limits_{k=0}^{j}(\bbbb)_kC_{k,r}(-j)_{k}(j-k)!=\frac{j!}{(1-m)_j}\sum\limits_{k=0}^{j}(-1)^k(\bbbb)_kC_{k,r}.
$$
Substituting these values into $q=m-1$ case of (\ref{eq:KRPTh1limit}) we arrive at (\ref{eq:KRPTh1limitm-1}).
Identity (\ref{eq:KRPThKummerq-m-1}) can be proved in a similar fashion or by applying confluence to (\ref{eq:KRPTh1limitm-1}).
$\hfill\square$

\textbf{Remark.}  In our concurrent  paper \cite[Theorem~5]{KPNew} we show that for $\cccc=\bbbb+1$
$$
Q^0_m(-j)=\frac{(\bbbb+1)_j}{(1-m)_j}{}_{r+2}F_{r+1}\!\left(\begin{matrix}-j,\bbbb,\f+\m\\\bbbb+1,\f\end{matrix}\right),
$$
so that (\ref{eq:KRPTh1limitm-1}) can be further simplified to
\begin{multline*}
(1-x)^{\aaaa}{}_{r+2}F_{r+1}\left.\!\left(\!\begin{matrix}\aaaa,\bbbb,\f+\m\\\bbbb+1, \f\end{matrix}\right\vert x\!\right)=
\sum\limits_{j=0}^{m-1}\frac{(\aaaa)_j}{j!}{}_{r+2}F_{r+1}\!\left(\begin{matrix}-j,\bbbb,\f+\m\\\bbbb+1,\f\end{matrix}\right)
\left(\frac{x}{x-1}\right)^j
\\
+\frac{x^{m}(\aaaa)_{m}(\f-\bbbb)_{\m}}{(x-1)^{m}(\bbbb+1)_{m}(\f)_{\m}}
{}_{2}F_{1}\left.\!\!\left(\begin{matrix}1, \aaaa+m\\\bbbb+m+1 \end{matrix}\right\vert\frac{x}{x-1}\!\right)
\end{multline*}
and similarly for (\ref{eq:KRPThKummerq-m-1}).

\begin{corollary}\label{cr:KRPTh3limitm-1}
Suppose $\aaaa+1\notin-\N_0$, $\aaaa-\bbbb\notin\{0,\ldots,m-1\}$ and $f_j-\aaaa\notin\{1-m_j,2-m_j,\ldots,0\}$ \emph{(}$1\le{j}\le{r}$\emph{)}. 
Then if $x\in\C\!\setminus\![1,\infty)$ the following identity holds\emph{:}
\begin{multline}\label{eq:KRPTh3limitm-1}
(1-x)^{\bbbb+m-1}{}_{r+2}F_{r+1}\left.\!\left(\!\begin{matrix}\aaaa,\bbbb,\f+\m\\\aaaa+1, \f\end{matrix}\right\vert x\!\right)=
\sum\limits_{j=0}^{m-1}\frac{(1-m)_j(\aaaa-\bbbb-m+1)_{j}}{(\aaaa+1)_{j}j!}\hat{Q}^0_m(-j)x^j
\\
+\frac{x^{m}(-1)^m(\bbbb)_m(\f-\aaaa)_{\m}}{(\aaaa+1)_{m}(\f)_{\m}}
{}_{2}F_{1}\left.\!\!\left(\begin{matrix}1, 1-\bbbb+\aaaa\\\aaaa+m+1 \end{matrix}\right\vert x\right),
\end{multline}
where
$$
\hat{Q}^0_m(-j)=\sum\limits_{k=0}^{j}\frac{(-1)^k(-j)_{k}(\aaaa)_k(\bbbb)_kC_{k,r}}{(1-m)_k(\aaaa-\bbbb-m+1)_k}
{}_{3}F_{2}\!\left(\begin{matrix}k-m, k-j,1-\bbbb+m\\k-m+1,\aaaa-\bbbb-m+k+1 \end{matrix}\!\right)
$$
and the numbers $C_{k,r}$ are defined in \emph{(\ref{eq:Ckr})}.
\end{corollary}
\textbf{Proof.}  For $q=m-1$ the polynomial $\hat{R}_{m-q-1}=\hat{R}_0$ defined in (\ref{eq:hatR}) reduces to a constant, namely:
\begin{multline*}
\hat{R}_0=\frac{(\bbbb)_{m}}{(\bbbb-\aaaa)_{m}}\sum\limits_{k=0}^{m-1}(-1)^{m-k}(\aaaa)_{k}C_{k,r}
+\frac{(\aaaa)_{m}(\bbbb)_{m}C_{m,r}}{(\bbbb-\aaaa)_{m}}
\\
=\frac{(-1)^{m}(\bbbb)_{m}}{(\bbbb-\aaaa)_{m}}\sum\limits_{k=0}^{m}(-1)^{k}(\aaaa)_{k}C_{k,r}
=\frac{(-1)^{m}(\bbbb)_{m}(\f-\aaaa)_{\m}}{(\bbbb-\aaaa)_{m}(\f)_{\m}},
\end{multline*}
where we used \cite[Corollary~3]{MP2012MathComm} in the last equality. The values $\hat{Q}^0_m(-j)$ are given in (\ref{eq:hatQm0negative}).
Substituting these values into $q=m-1$ case of (\ref{eq:KRPTh2limit}) we arrive at (\ref{eq:KRPTh3limitm-1}). $\hfill\square$

\textbf{Remark.}   Both formulas (\ref{eq:KRPTh1limitm-1}) and (\ref{eq:KRPTh3limitm-1}) can be viewed as extensions of Karlsson's summation theorem \cite[(4.2)]{MS2010}.
Indeed,  Karlsson's formula  can be recovered from (\ref{eq:KRPTh1limitm-1}) by setting $x=1$ and using the well known asymptotic approximation for ${}_2F_{1}$ on the right hand side. Similarly, setting $x=1$ in (\ref{eq:KRPTh1limitm-1}) and applying Euler's transformation and the Gauss summation formula on the right hand side we again arrive at Karlsson's result.

Finally,  we treat the case $q=m-2$ when the polynomials $R_{m-q-1}$ and $\hat{R}_{m-q-1}$ become linear, so that their roots can be explicitly calculated.

\begin{corollary}\label{cr:KRPTh1limitm-2}
If $\bbbb+2\notin-\N_0$ and $f_j-\bbbb\notin\{2-m_j,3-m_j,\ldots,0\}$ for indices $0\le{j}\le{r}$ such that $m_j\ge2$, then the following identity holds\emph{:}
\begin{multline}\label{eq:KRPTh1limitm-2}
(1-x)^{\aaaa}{}_{r+2}F_{r+1}\left.\!\left(\begin{matrix}\aaaa,\bbbb,\f+\m\\\bbbb+2, \f\end{matrix}\right\vert x\right)=
\sum\limits_{j=0}^{m-2}\frac{(\aaaa)_j}{(\bbbb+2)_j}\left(\frac{x}{x-1}\right)^j\sum\limits_{k=0}^{j}(-1)^k(\bbbb)_k(j-k+1)C_{k,r}
\\
+\frac{x^{m-1}(\aaaa)_{m-1}\left((\bbbb+m)(\f-\bbbb)_{\m}-\bbbb(\f-\bbbb-1)_{\m}\right)}{(x-1)^{m-1}(\bbbb+2)_{m-1}(\f)_{\m}}
{}_{3}F_{2}\left.\!\!\left(\begin{matrix}1,\aaaa+m-1,\lambda+m\\\bbbb+m+1, \lambda+m-1\end{matrix}\right\vert\frac{x}{x-1}\!\right)
\end{multline}
for all $x\in\C\!\setminus\![1,\infty)$ and
\begin{multline}\label{eq:KRPThKummerq-m-2}
e^{-x}{}_{r+1}F_{r+1}\left.\left(\begin{matrix}\bbbb,\f+\m\\\bbbb+2, \f\end{matrix}\right\vert x\right)=
\sum\limits_{j=0}^{m-2}\frac{(-x)^j}{(\bbbb+2)_j}\sum\limits_{k=0}^{j}(-1)^k(\bbbb)_k(j-k+1)C_{k,r}
\\
+\frac{(-x)^{m-1}\left((\bbbb+m)(\f-\bbbb)_{\m}-\bbbb(\f-\bbbb-1)_{\m}\right)}{(\bbbb+2)_{m-1}(\f)_{\m}}
{}_{2}F_{2}\left.\!\left(\begin{matrix}1,\lambda+m\\\bbbb+m+1, \lambda+m-1\end{matrix}\right\vert-x\right)
\end{multline}
for all $x\in\C$.  Here
\begin{equation}\label{eq:R1root}
\lambda=1+\bbbb-\frac{\bbbb(\f-\bbbb-1)}{(\f-\bbbb-1+\m)}
\end{equation}
and the numbers $C_{k,r}$ are defined in \emph{(\ref{eq:Ckr})}.
\end{corollary}
\textbf{Proof.} Indeed, taking $q=m-2$ in (\ref{eq:KRPTh1limit}) we get formula (\ref{eq:KRPTh1limitm-2})  with $Q_m^0(-j)$ in the first sum and the factor in front of ${}_3F_{2}$
given by
$$
\frac{x^{m-1}(\aaaa)_{m-1}R_{1}(1-m)}{(x-1)^{m-1}(\bbbb+2)_{m-1}}.
$$
The number  $\lambda$ is the root of the equation $R_{1}(t)=0$, where $R_{1}(t)$ is defined by $q=m-2$ case of (\ref{eq:R}).  Using the form (\ref{eq:Rnew}) for $R_{1}(t)$
we easily get formula (\ref{eq:R1root}) for the root and the value
$$
R_1(1-m)=\frac{(\bbbb+m)(f-\bbbb)_{\m}-\bbbb(\f-\bbbb-1)_{\m}}{(\f)_{\m}}.
$$
Further according to (\ref{eq:Qm0negative}) with  $q=m-2$ we have
$$
Q^0_m(-j)=\frac{1}{(2-m)_j}\sum\limits_{k=0}^{j}(\bbbb)_kC_{k,r}(-j)_k(2)_{j-k}=\frac{j!}{(2-m)_j}\sum\limits_{k=0}^{j}(-1)^k(\bbbb)_kC_{k,r}(j-k+1),
$$
which yields (\ref{eq:KRPTh1limitm-2}). Identity (\ref{eq:KRPThKummerq-m-2}) can be proved in a similar fashion or by applying confluence to (\ref{eq:KRPTh1limitm-2}).
$\hfill\square$

\textbf{Remark.}  In our concurrent paper \cite[Theorem~5]{KPNew} we show that for $\cccc=\bbbb+2$
$$
Q^0_m(-j)=\frac{(\bbbb+2)_j}{(2-m)_j}{}_{r+2}F_{r+1}\!\left(\begin{matrix}-j,\bbbb,\f+\m\\\bbbb+2,\f\end{matrix}\right)
$$
so that (\ref{eq:KRPTh1limitm-2}) can be further simplified to
\begin{multline*}
(1-x)^{\aaaa}{}_{r+2}F_{r+1}\left.\!\left(\!\begin{matrix}\aaaa,\bbbb,\f+\m\\\bbbb+2, \f\end{matrix}\right\vert x\!\right)=
\sum\limits_{j=0}^{m-2}\frac{(\aaaa)_j}{j!}\left(\frac{x}{x-1}\right)^j{}_{r+2}F_{r+1}\!\left(\begin{matrix}-j,\bbbb,\f+\m\\\bbbb+2,\f\end{matrix}\right)
\\
+\frac{x^{m-1}(\aaaa)_{m-1}\left[(\bbbb+m)(\f-\bbbb)_{\m}-\bbbb(\f-\bbbb-1)_{\m}\right]}{(x-1)^{m-1}(\bbbb+2)_{m-1}(\f)_{\m}}
{}_{3}F_{2}\left.\!\!\left(\begin{matrix}1,\aaaa+m-1,\lambda+m\\\bbbb+m+1, \lambda+m-1\end{matrix}\right\vert\frac{x}{x-1}\!\right)
\end{multline*}
and similarly for  (\ref{eq:KRPThKummerq-m-2}).

\begin{corollary}\label{cr:KRPTh2limitm-2}
Suppose $\aaaa+2\notin-\N_0$,  $\aaaa-\bbbb\notin\{-1,\ldots,m-2\}$ and $f_j-\aaaa\notin\{2-m_j,3-m_j,\ldots,0\}$ for indices $0\le{j}\le{r}$ such that $m_j\ge2$. 
Then for all $x\in\C\!\setminus\![1,\infty)$ the following identity holds\emph{:}
\begin{multline}\label{eq:KRPTh2limitm-2}
(1-x)^{\bbbb+m-2}{}_{r+2}F_{r+1}\left.\!\left(\!\begin{matrix}\aaaa,\bbbb,\f+\m\\\aaaa+2, \f\end{matrix}\right\vert x\!\right)=\sum\limits_{j=0}^{m-2}\frac{(2-m)_j(\aaaa-\bbbb-m+2)_j}{(\aaaa+2)_jj!}
\hat{Q}^0_m(-j)x^j
\\
+x^{m-1}\frac{(-1)^{m-1}(\bbbb)_{m-1}[(\aaaa+m)(\f-\aaaa)_{\m}-\aaaa(\f-\aaaa-1)_{\m}]}{(\aaaa+2)_{m-1}(\f)_{\m}}
{}_{3}F_{2}\left.\!\left(\!\begin{matrix}1,\aaaa-\bbbb+1,\gamma+m\\\aaaa+m+1, \gamma+m-1\end{matrix}\right\vert x\!\right),
\end{multline}
where
\begin{equation}\label{eq:Rhat1root}
\gamma=1-m+\frac{\aaaa(\f-\aaaa-1)-(\aaaa+m)(\f-\aaaa-1+\m)}{\aaaa(\f-\aaaa-1)/(\aaaa-\bbbb+1)-(\f-\aaaa-1+\m)}
\end{equation}
and the values $\hat{Q}^0_m(-j)$, $j=0,\ldots,m-2$ are given by \emph{(\ref{eq:hatQm0negative})} with $q=m-2$.
\end{corollary}
\textbf{Proof.}  Indeed, taking $q=m-2$ in (\ref{eq:KRPTh2limit}) we get formula (\ref{eq:KRPTh2limitm-2}) with the factor in front of ${}_3F_{2}$
given by
$$
x^{m-1}\frac{\hat{R}_{1}(1-m)(\bbbb-\aaaa)_{m-1}}{(\aaaa+2)_{m-1}}
$$
and $\gamma$ being the root of $\hat{R}_{1}(t)=0$, where $\hat{R}_{1}(t)$ is defined by $q=m-2$ case of (\ref{eq:hatR}), i.e.
\begin{multline*}
\hat{R}_{1}(t)=\frac{(\bbbb)_{m-1}}{(\bbbb-\aaaa)_{m-1}}\sum\limits_{k=0}^{m-2}\frac{(-m+k)_{m-1-k}(\aaaa)_{k}C_{k,r}}{(m-k-1)!}
{}_{3}F_{2}\!\!\left(\begin{matrix}-1,t+m-1,1-\bbbb-k\\1-\bbbb+\aaaa,m-k\end{matrix}\right)
\\
+\sum\limits_{k=m-1}^{m}\frac{(-1)^k(\aaaa)_{k}(\bbbb)_{k}C_{k,r}(t+m-1)_{k-m+1}}{(\aaaa-\bbbb-m+2)_k(k-m+1)!}
{}_{3}F_{2}\!\!\left(\begin{matrix}-m+k,t+k,-\bbbb-m+2\\\aaaa-\bbbb-m+2+k,k-m+2\end{matrix}\right).
\end{multline*}
Collecting the terms containing $(t+m-1)$ and using $(-m+k)_{m-1-k}=(-1)^{m-1-k}(m-k)!$ we get:
\begin{multline*}
\hat{R}_{1}(t)=\frac{(\bbbb)_{m-1}(-1)^{m-1}}{(\bbbb-\aaaa)_{m-1}}\sum\limits_{k=0}^{m-1}(-1)^{k}(\aaaa)_{k}C_{k,r}(m-k)
\\
+\frac{(t+m-1)(\bbbb)_{m-1}(-1)^{m-1}}{(\bbbb-\aaaa-1)_{m}}\sum\limits_{k=0}^{m}(\aaaa)_{k}(-1)^{k}C_{k,r}(1-\bbbb-k)
\\
=\frac{(\bbbb)_{m-1}(-1)^{m-1}}{(\bbbb-\aaaa)_{m-1}}R_1(1-m)+\frac{(t+m-1)(\bbbb)_{m-1}(-1)^{m-1}}{(\bbbb-\aaaa-1)_{m}}R_1(\bbbb),
\end{multline*}
where $R_1(t)$ is the polynomial defined in (\ref{eq:R}) (with $\aaaa$ in place of $\bbbb$).  Using the form (\ref{eq:Rnew})
for $R_{1}(t)$ we can solve $\hat{R}_{1}(t)=0$ to get (\ref{eq:Rhat1root}). According to (\ref{eq:hatR-q-1}) we have
\begin{multline*}
\hat{R}_{1}(1-m)=\frac{(\bbbb)_{m-1}}{(\bbbb-\aaaa)_{m-1}}\sum\limits_{k=0}^{m-1}\frac{(-m+k)_{m-1-k}(\aaaa)_kC_{k,r}}{(m-1-k)!}
\\
=\frac{(\bbbb)_{m-1}(-1)^{m-1}}{(\bbbb-\aaaa)_{m-1}}\sum\limits_{k=0}^{m-1}(-1)^{k}(m-k)(\aaaa)_kC_{k,r}=\frac{(\bbbb)_{m-1}(-1)^{m-1}}{(\bbbb-\aaaa)_{m-1}}R_1(1-m)
\\
=\frac{(\bbbb)_{m-1}(-1)^{m-1}[(\bbbb+m)(\f-\bbbb)_{\m}-\bbbb(\f-\bbbb-1)_{\m}]}{(\bbbb-\aaaa)_{m-1}(\f)_{\m}},
\end{multline*}
where we used the last formula in the proof of Corollary~\ref{cr:KRPTh1limitm-2}.$\hfill\square$

\section{Summation and transformation formulas for $x=1$}

The limit case of Theorem~\ref{th:KRPTh1limit} leads to the next generalization of Karlsson's summation theorem \cite[(4.2)]{MS2010}.
\begin{theorem}\label{th:genKarlsson}
Suppose $q\in\{0,\ldots,m-1\}$ and $\Re(\aaaa+q)<0$. Then the next summation formula holds
\begin{equation}\label{eq:genKarlsson}
{}_{r+2}F_{r+1}\!\!\left(\begin{matrix}\aaaa,\bbbb,\f+\m\\\bbbb+m-q, \f\end{matrix}\right)
=\frac{\Gamma(\bbbb+m-q)\Gamma(1-\aaaa)R_{m-q-1}(\aaaa)}{\Gamma(\bbbb-\aaaa+m-q)(m-q-1)!},
\end{equation}
where $R_{m-q-1}(t)$ is given by  \emph{(\ref{eq:R})} or \emph{(\ref{eq:Rnew})}.
\end{theorem}

\noindent\textbf{Remark.} Formula (\ref{eq:Rnew}) shows, in particular, that the right hand side vanishes if $m-q-m_j\le0$ for some $0\le{j}\le{r}$ and $f_j-\bbbb\in\{m-q-m_j,m-q-m_j+1,\ldots,0\}$. Hence, so does the left hand side.  We further remark that we give a different proof and generalization of 
(\ref{eq:genKarlsson}) in our concurrent paper \cite[Theorem~2]{KPNew}.

\noindent\textbf{Proof.} We will need the asymptotic formula \cite[(16.11.6)]{NIST}
\begin{equation}\label{eq:p+1Fp-asymp}
\frac{\Gamma(\a)}{\Gamma(\b)}{}_{p+1}F_{p}\left.\!\left(\begin{matrix}\a\\\b\end{matrix}\right\vert -z\right)=\sum_{k=1}^{p+1}\frac{\Gamma(a_k)\Gamma(\a_{[k]}-a_k)}{\Gamma(\b-a_k)}z^{-a_k}\bigl(1+\O(1/z)\bigr),~~\text{as}~|z|\to\infty,
\end{equation}
valid if $|\arg(z)|<\pi$ and $a_j-a_k\notin\Z$ for all $k\ne{j}$.  Denote $-z=x/(x-1)$, assume that $\Re(\aaaa+q)<0$ and the upper parameters of ${}_{m-q+1}F_{m-q}$
on the right hand side of (\ref{eq:KRPTh1limit}) do not differ by integers.  Then, after multiplication by $(1-x)^{-\aaaa}$, formula (\ref{eq:KRPTh1limit}) takes the form:
\begin{multline*}
{}_{r+2}F_{r+1}\left.\!\left(\begin{matrix}\aaaa,\bbbb,\f+\m\\\bbbb+m-q, \f\end{matrix}\right\vert x\!\right)=
\sum\limits_{j=0}^{q}\frac{(\aaaa)_j(-q)_jQ^{0}_m(-j)}{(\bbbb+m-q)_jj!}(-z)^j(1+z)^{\aaaa}
\\
+z^{q+1}(1+z)^{\aaaa}\frac{(-1)^{q+1}(\aaaa)_{q+1}R_{m-q-1}(-q-1)}{(\bbbb+m-q)_{q+1}(m-q-1)!}
{}_{m-q+1}F_{m-q}\left.\!\!\left(\begin{matrix}1,\aaaa+q+1,\llambda+q+2\\\bbbb+m+1, \llambda+q+1\end{matrix}\right\vert-z\!\right).
\end{multline*}
Clearly $z\to+\infty$ as $x\to1-$. Under the above assumptions the first sum on the right vanishes as $z\to+\infty$ and, according to the asymptotic formula (\ref{eq:p+1Fp-asymp}), \begin{multline*}
z^{q+1}(1+z)^{\aaaa}{}_{m-q+1}F_{m-q}\left.\!\!\left(\begin{matrix}1,\aaaa+q+1,\llambda+q+2\\\bbbb+m+1, \llambda+q+1\end{matrix}\right\vert-z\!\right)
\\
=\left(C_1+C_2z^{\aaaa+q}+C_3z^{\aaaa-\lambda_1-1}+\cdots+C_{m-q+1}z^{\aaaa-\lambda_{m-q-1}-1}\right)\bigl(1+\O(1/z)\bigr),
\end{multline*}
where the constants $C_1,\ldots,C_{m-q+1}$ can be read off (\ref{eq:p+1Fp-asymp}).  Assuming further that $\Re(\aaaa-\lambda_k-1)<0$
for $k=1,\ldots,m-q-1$, we conclude that the limit of the above expression is equal to
\begin{multline*}
C_1=\frac{\Gamma(\bbbb+m+1)\Gamma(\llambda+q+1)\Gamma(-\aaaa-q)\Gamma(\llambda-\aaaa+1)}{\Gamma(\llambda+q+2)\Gamma(\bbbb-\aaaa+m-q)\Gamma(\llambda-\aaaa)}
\\
=\frac{\Gamma(\bbbb+m+1)\Gamma(-\aaaa-q)(\llambda-\aaaa)_1}{\Gamma(\bbbb-\aaaa+m-q)(\llambda+q+1)_1}
=\frac{\Gamma(\bbbb+m+1)\Gamma(-\aaaa-q)R_{m-q-1}(\aaaa)}{\Gamma(\bbbb-\aaaa+m-q)R_{m-q-1}(-q-1)},
\end{multline*}
which leads to (\ref{eq:genKarlsson}) once we apply $\Gamma(-\aaaa-q)=(-1)^{q+1}\Gamma(1-\aaaa)/(\aaaa)_{q+1}$.  The auxiliary assumptions $\Re(\aaaa-\lambda_k-1)<0$
  can now be removed by resorting to analytic continuation.  $\hfill\square$

\smallskip

In \cite{KRP2014} the authors applied the beta integral method to formulas (\ref{eq:KRPTh1-1}) and (\ref{eq:KRPTh1-2}) to obtain transformation formulas
for terminating and non-terminating ${}_{p+1}F_{p}(1)$.  We can similarly apply the beta integral method to  (\ref{eq:KRPTh1limit})
and (\ref{eq:KRPTh2limit}) to get the transformation formulas presented below.

\begin{theorem}\label{th:Thomaelike1}
Suppose $q\in\{0,\ldots,m-1\}$, $\bbbb+m-q\notin-\N_0$ and 
$f_j-\bbbb\notin\{m-q-m_j,m-q-m_j+1,\ldots,0\}$ for indices $0\le{j}\le{r}$ such that $m-q-m_j\le0$. Then for each $n\in\N$ the next transformation formula holds\emph{:}
\begin{multline}\label{eq:Thomaelike1}
{}_{r+3}F_{r+2}\!\!\left(\begin{matrix}-n,\bbbb,\dddd,\f+\m\\\bbbb+m-q, \eeee, \f\end{matrix}\right)
\!=\!\frac{(\eeee-\dddd)_n}{(\eeee)_n}\!\sum\limits_{j=0}^{q}\!\frac{(-n)_j(-q)_j(\dddd)_jQ^0_m(-j)}{(\bbbb+m-q)_{j}(1+\dddd-\eeee-n)_{j}j!}
+\frac{(-1)^{q+1}(-n)_{q+1}}{(\bbbb+m-q)_{q+1}}
\\
\times\!\frac{(\eeee-\dddd)_{n-q-1}R_{m-q-1}(-q-1)}{(\eeee)_{n}(m-q-1)!}{}_{m-q+2}F_{m-q+1}\!\left(\begin{matrix}-n+q+1,1,\dddd+q+1,\llambda+q+2\\\bbbb+m+1,2+\dddd+q-\eeee-n, \llambda+q+1\end{matrix}\!\right),
\end{multline}
where $\llambda$ are the roots of $R_{m-q-1}(t)$ given by \emph{(\ref{eq:R})} or \emph{(\ref{eq:Rnew})},  the values $Q^0_m(-j)$ are given in \emph{(\ref{eq:Qm0negative})} and $R_{m-q-1}(-q-1)$
is computed in \emph{(\ref{eq:R-q-1})}. If $n\le{q}$ the second term is missing.
\end{theorem}
\textbf{Proof.}  The proof is by straightforward application of the beta integral method: set $\aaaa=-n$, multiply both sides of (\ref{eq:KRPTh1limit}) by $x^{\dddd-1}(1-x)^{n+\eeee-\dddd-1}$
and integrate with respect to $x$ over the interval $[0,1]$. $\hfill\square$

\smallskip

\begin{theorem}\label{th:Thomaelike2}
Suppose $q\in\{0,\ldots,m-1\}$, $\aaaa+m-q\notin-\N_0$, $\aaaa-\bbbb\notin\{q+1-m,\ldots,q\}$ and 
$f_j-\aaaa\notin\{m-q-m_j,m-q-m_j+1,\ldots,0\}$ for indices $0\le{j}\le{r}$ such that $m-q-m_j\le0$.
Then the next transformation formula holds\emph{:}
\begin{multline}\label{eq:Thomaelike2}
{}_{r+3}F_{r+2}\!\!\left(\begin{matrix}\aaaa,\bbbb,\dddd,\f+\m\\\aaaa+m-q, \eeee, \f\end{matrix}\right)
\!=\!\frac{\Gamma(\eeee)\Gamma(\eeee-\bbbb-\dddd-q)}{\Gamma(\eeee-\dddd)\Gamma(\eeee-\bbbb-q)}\!\sum\limits_{j=0}^{q}\frac{(-q)_j(\aaaa-\bbbb-q)_j(\dddd)_j\hat{Q}^0_m(-j)}{(\aaaa+m-q)_{j}(\eeee-\bbbb-q)_{j}j!}
\\
+\frac{\Gamma(\eeee)\Gamma(\eeee-\bbbb-\dddd-q)(\bbbb-\aaaa)_{q+1}(\dddd)_{q+1}\hat{R}_{m-q-1}(-q-1)}{\Gamma(\eeee-\dddd)\Gamma(\eeee-\bbbb+1)(\aaaa+m-q)_{q+1}}
\\
\times{}_{m-q+2}F_{m-q+1}\!\left(\begin{matrix}1,\aaaa-\bbbb+1, \dddd+q+1, \ggamma+q+2\\\aaaa+m+1, \eeee-\bbbb+1, \ggamma+q+1\end{matrix}\!\right),
\end{multline}
where $\ggamma$ are the roots of $\hat{R}_{m-q-1}(t)$ given by \emph{(\ref{eq:hatR})}, the values $\hat{Q}^0_m(-j)$ are given in \emph{(\ref{eq:hatQm0negative})} and $\hat{R}_{m-q-1}(-q-1)$ is computed in \emph{(\ref{eq:hatR-q-1})}.
\end{theorem}
\textbf{Proof.}  The proof is by straightforward application of the beta integral method: multiply both sides of (\ref{eq:KRPTh2limit}) by $x^{\dddd-1}(1-x)^{\eeee-\dddd-\bbbb-q-1}$
and integrate with respect to $x$ over the interval $[0,1]$ termwise. $\hfill\square$

\smallskip

The cases $q=0$ of the above two theorems are given in the next corollary.
\begin{corollary}\label{cr:Thomae-q0}
For $n\in\N$, $\bbbb+m\notin-\N_0$ and $f\ne{\bbbb}$ if $m=r=1$, we have
\begin{equation}\label{eq:eq:KRP1limit}
\frac{(\eeee)_n}{(\eeee-\dddd)_n}{}_{r+3}F_{r+2}\!\!\left(\!\begin{matrix}-n, \bbbb, \dddd, \f+\m\\\bbbb+m, \eeee, \f\end{matrix}\!\right)
\!\!=\!\!\frac{(\bbbb)_m}{(\f)_{\m}}+\left(1-\frac{(\bbbb)_m}{(\f)_{\m}}\right){}_{m+2}F_{m+1}\!\!\left(\begin{matrix}-n, 1, \dddd, \llambda+1\\\bbbb+m, 1-\eeee+\dddd-n, \llambda\end{matrix}\right),
\end{equation}
where $\llambda$ are the roots of the polynomial $R_{m-1}(t)$ given in \emph{(\ref{eq:Rm-1})}. For $\aaaa+m\notin-\N_0$, $\aaaa-\bbbb\notin\{1-m,\ldots,0\}$ and $f\ne{\aaaa}$ if $m=r=1$, we have
\begin{multline}\label{eq:KRP2limit}
\frac{\Gamma(\eeee-\dddd)}{\Gamma(\eeee-\bbbb-\dddd)\Gamma(\dddd)}{}_{r+3}F_{r+2}\!\left(\begin{matrix}\aaaa,\bbbb,\dddd,\f+\m\\\aaaa+m,\eeee,\f\end{matrix}\right)
\\
=1+\bbbb \dddd\biggl(\frac{\aaaa(\f+\m)}{(\aaaa+m)(\f)}-1\biggr)
{}_{m+2}F_{m+1}\!\left(\begin{matrix}1, \aaaa-\bbbb+1, \dddd+1, \ggamma+2\\\aaaa+m+1,\eeee-\bbbb+1, \ggamma+1\end{matrix}\right),
\end{multline}
where $\ggamma$ are the roots of the polynomial $\hat{R}_{m-1}(t)$ given in \emph{(\ref{eq:hatRm-1})}.
\end{corollary}

\medskip

We conclude the paper with some examples of transformation and summation formulas established above.

\textbf{Example~1.} For $r=1$, the numbers $C_{m,r}$ are $C_{0,1}=1$, $C_{1,1}=2/f$, $C_{2,1}=1/(f)_2$, so that
$$
R_1(t)=(1-t)+(2\bbbb/f)t-(\bbbb)_2/(f)_2(t+1)
$$
and
$$
\lambda_1=\frac{f+\bbbb+1}{f-\bbbb+1}.
$$
Formula (\ref{eq:KRPTh1-1-reduced}) then reduces to
$$
(f)_2(1-x)^{\aaaa}{}_{3}F_{2}\left.\!\left(\begin{matrix}\aaaa, \bbbb, f+2\\\bbbb+2,f\end{matrix}\right\vert x\right)
=(\bbbb)_2+[(f)_2-(\bbbb)_2]{}_{3}F_{2}\left.\!\!\!\left(\begin{matrix}\aaaa, 1, \lambda_1+1\\\bbbb+2, \lambda_1\end{matrix}\right\vert\frac{x}{x-1}\right),
$$
while formula (\ref{eq:eq:KRP1limit}) takes the form
$$
\frac{(\eeee)_n(f)_2}{(\eeee-\dddd)_n}{}_{4}F_{3}\!\left(\!\begin{matrix}-n, \bbbb, \dddd, f+2\\\bbbb+2, \eeee, f\end{matrix}\!\right)
\!\!=\!\!(\bbbb)_2+[(f)_2-(\bbbb)_2]{}_{4}F_{3}\!\!\left(\!\begin{matrix}-n, 1, \dddd, \lambda_1+1\\\bbbb+2, 1-\eeee+\dddd-n, \lambda_1\end{matrix}\!\right).
$$

\medskip

\textbf{Example~2.} This example illustrates the remark made in the introduction around formula (\ref{eq:limit-m1}).
Start with \cite[p.116]{KRP2014}
$$
{}_{4}F_{3}\!\left(\begin{matrix}\aaaa,\bbbb, \dddd, f+1\\\cccc, \eeee, f\end{matrix}\right)
=\frac{\Gamma(\eeee)\Gamma(\psi)}{\Gamma(\eeee-\dddd)\Gamma(\psi+\dddd)}
{}_{4}F_{3}\!\left(\begin{matrix}\dddd, \cccc-\aaaa-1, \cccc-\bbbb-1, \eta+1\\\cccc,\psi+\dddd, \eta\end{matrix}\right)
$$
where $\psi=\cccc+\eeee-\aaaa-\bbbb-\dddd-m$ and
$$
\eta=\frac{(\cccc-\aaaa-1)(\cccc-\bbbb-1)f}{\aaaa \bbbb+(\cccc-\aaaa-\bbbb-1)f}.
$$
Letting $\cccc-\aaaa-1\to0$ in the above formula and applying (\ref{eq:limit-m1}) with $x=1$ we get
$$
{}_{4}F_{3}\!\left(\begin{matrix}\aaaa,\bbbb, \dddd, f+1\\\aaaa+1,\eeee,f\end{matrix}\right)
=\frac{\Gamma(\eeee)\Gamma(\eeee-\bbbb-\dddd)}{\Gamma(\eeee-\dddd)\Gamma(\eeee-\bbbb)}
\left[1-\frac{\bbbb(\aaaa-f)}{(\aaaa-\bbbb)f}\left(1-{}_{3}F_{2}\!\left(\begin{matrix}\dddd, \aaaa-\bbbb, 1\\\aaaa+1,\eeee-\bbbb\end{matrix}\right)\right)\right].
$$
For $\eeee=\bbbb+1$ this leads to known summation formula \cite[7.5.3.1]{PBM3}.  For $\eeee=\bbbb+r$, with integer $r\ge1$, we will apply the easily verifiable identity
\cite[Lemma~4]{MP2012}
\begin{equation}\label{eq:pFq1-r}
{}_{p+1}F_{q+1}\!\left(\begin{matrix}\a,1\\\b,r\end{matrix}\right)=\frac{(1-\b)_{r-1}(r-1)!}{(1-\a)_{r-1}(-1)^{(r-1)(p-q)}}\left[{}_{p}F_{q}\!\left(\begin{matrix}\a-r+1\\\b-r+1\end{matrix}\right)
-\sum_{j=0}^{r-2}\frac{(\a-r+1)_j}{(\b-r+1)_jj!}\right]
\end{equation}
to get by invoking the Chu-Vandermonde identity \cite[Corollary~2.2.3]{AAR} a presumably new summation formula:
\begin{multline*}
{}_{4}F_{3}\!\left(\begin{matrix}\aaaa,\bbbb,\dddd,f+1\\\aaaa+1,\bbbb+r,f\end{matrix}\right)
=\frac{\Gamma(\eeee)\Gamma(\eeee-\bbbb-\dddd)}{\Gamma(\eeee-\dddd)\Gamma(r)}
\left[1-\frac{\bbbb(\aaaa-f)}{(\aaaa-\bbbb)f}\left(1-{}_{3}F_{2}\!\left(\begin{matrix}\dddd, \aaaa-\bbbb, 1\\\aaaa+1,r\end{matrix}\right)\right)\right]
\\
=\frac{\Gamma(\bbbb+r)\Gamma(r-\dddd)}{\Gamma(\bbbb-\dddd+r)\Gamma(r)}
\biggl[1-\frac{\bbbb(\aaaa-f)}{(\aaaa-\bbbb)f}\biggl(1-\frac{(-1)^{r-1}(-\aaaa)_{r-1}(r-1)!}{(1-\dddd)_{r-1}(1-\aaaa+\bbbb)_{r-1}}\times
\\
\biggl\{\frac{\Gamma(\aaaa-r+2)\Gamma(\bbbb-\dddd+r)}{\Gamma(1-\dddd+\aaaa)\Gamma(1+\bbbb)}
-\sum_{j=0}^{r-2}\frac{(\dddd-r+1)_j(\aaaa-\bbbb-r+1)_j}{(\aaaa-r+2)_jj!}\biggr\}\biggr)\biggr].
\end{multline*}

\medskip

\textbf{Example~3.} For $r=1$, $m=2$ using $C_{k,1}=\binom{m}{k}/(f)_k$ we have:
$$
\hat{R}_1(t)=\frac{2\bbbb}{\aaaa-\bbbb}\left(1-\frac{(t+1)(1-\bbbb)}{2(\aaaa-\bbbb+1)}\right)-\frac{2\aaaa \bbbb}{f(\aaaa-\bbbb)}\left(1+\frac{(t+1)\bbbb}{\aaaa-\bbbb+1}\right)+\frac{(\aaaa)_2(\bbbb)_2(t+1)}{(f)_2(\aaaa-\bbbb)_2}.
$$
We can solve $\hat{R}_1(t)=0$ or simplify the expression for the root calculated in (\ref{eq:Rhat1root}) to get:
$$
\gamma=\frac{(1+\aaaa-\bbbb)(1+f)+\aaaa(f-\bbbb)}{(1+\aaaa-\bbbb)(1+f)-\aaaa(f-\bbbb)}.
$$
Formula (\ref{eq:KRPTh1-2-reduced}) takes the form:
$$
(1-x)^{\bbbb}{}_{3}F_{2}\left.\!\left(\!\begin{matrix}\aaaa,\bbbb,f+2\\\aaaa+2, f\end{matrix}\right\vert x\!\right)
=1
+x\bbbb\biggl(\frac{\aaaa(f+2)}{(\aaaa+2)f}-1\biggr)
{}_{3}F_{2}\left.\!\left(\!\begin{matrix}1,\aaaa-\bbbb+1,\gamma+2\\\aaaa+3, \gamma+1\end{matrix}\right\vert x\!\right),
$$
where $\aaaa-\bbbb\notin\{-1,0\}$.  Note that the alternative form (\ref{eq:KRPTh1-2-reducedAlt}) of this formula is:
$$
(1-x)^{\bbbb}{}_{3}F_{2}\left.\!\!\left(\begin{matrix}\aaaa,\bbbb,f+2\\\aaaa+2, f\end{matrix}\right\vert x\right)
=A+(1-A){}_{3}F_{2}\left.\!\!\left(\begin{matrix}\aaaa-\bbbb, 1, \gamma+1\\\aaaa+2, \gamma\end{matrix}\right\vert x\right)
$$
with
$$
A=\sum\limits_{k=0}^{2}\binom{2}{k}\frac{(-1)^k(\aaaa)_k(\bbbb)_k(-1-\aaaa)_{2-k}}{(f)_k(\aaaa-\bbbb)_k(\bbbb-\aaaa-1)_{2-k}}.
$$
Formula (\ref{eq:KRP2limit}) takes the form:
$$
{}_{4}F_{3}\!\left(\begin{matrix}\aaaa,\bbbb,\dddd,f+2\\\aaaa+2,\eeee,f\end{matrix}\right)
\!=\!\frac{\Gamma(\dddd)\Gamma(\eeee-\bbbb-\dddd)}{\Gamma(\eeee-\dddd)}
\!\left[1+\bbbb \dddd\biggl(\frac{\aaaa(f+2)}{(\aaaa+2)f}-1\!\biggr)
{}_{4}F_{3}\!\left(\begin{matrix}1, \aaaa-\bbbb+1, \dddd+1, \gamma+2\\\aaaa+3,\eeee-\bbbb+1, \gamma+1\end{matrix}\right)\right]
$$
or
$$
{}_{4}F_{3}\!\left(\begin{matrix}\aaaa,\bbbb,\dddd,f+2\\\aaaa+2,\eeee,f\end{matrix}\right)
=\frac{\Gamma(\eeee)\Gamma(\eeee-\bbbb-\dddd)}{\Gamma(\eeee-\dddd)\Gamma(\eeee-\bbbb)}\biggl[A+(1-A){}_{4}F_{3}\!\left(\begin{matrix}\dddd, \aaaa-\bbbb, 1, \gamma+1\\\aaaa+2,\eeee-\bbbb, \gamma\end{matrix}\right)\biggr].
$$
In particular, when $\eeee=\bbbb+1$ the last transformation leads to the summation formula:
\begin{multline*}
{}_{4}F_{3}\!\left(\begin{matrix}\aaaa,\bbbb,\dddd,f+2\\\aaaa+2,\bbbb+1,f\end{matrix}\right)
=\frac{\Gamma(\bbbb+1)\Gamma(1-\dddd)}{\Gamma(\eeee-\dddd)}\biggl[A+(1-A){}_{3}F_{2}\!\left(\begin{matrix}\dddd, \aaaa-\bbbb, \lambda+1\\\aaaa+2, \lambda\end{matrix}\right)\biggr]
\\
=\frac{\Gamma(\bbbb+1)\Gamma(1-\dddd)}{\Gamma(\eeee-\dddd)}\biggl[A+(1-A)
\frac{\Gamma(\aaaa+2)\Gamma(\bbbb-\dddd+2)}{\Gamma(\bbbb+2)\Gamma(\aaaa-\dddd+2)}\biggl(1-\frac{\dddd(\aaaa-\bbbb)}{\lambda(\dddd-\bbbb-1)}\biggr)\biggr],
\end{multline*}
where we applied \cite[7.4.4.9]{PBM3} (see also \cite[Lemma~1]{Miller2005}).  If $\eeee=\bbbb+2$ by using the $r=2$ case of (\ref{eq:pFq1-r}) we obtain:
\begin{multline*}
{}_{4}F_{3}\!\left(\begin{matrix}\aaaa,\bbbb,\dddd,f+2\\\aaaa+2,\bbbb+2,f\end{matrix}\right)
=\frac{\Gamma(\bbbb+2)\Gamma(2-\dddd)}{\Gamma(\bbbb-\dddd+2)}\times
\\
\biggl[A+\frac{(1-A)(\aaaa+1)(\lambda-1)}{\lambda(\dddd-1)(\aaaa-\bbbb-1)}\left(
\frac{\Gamma(\aaaa+1)\Gamma(\bbbb-\dddd+3)}{\Gamma(\bbbb+2)\Gamma(\aaaa-\dddd+2)}\left(1-\frac{(\dddd-1)(\aaaa-\bbbb-1)}{(\lambda-1)(\dddd-\bbbb-2)}\right)-1\right)\biggr].
\end{multline*}
In a similar fashion setting $\eeee=\bbbb+r$, $r>2$, and applying (\ref{eq:pFq1-r}) with \cite[7.4.4.9]{PBM3} we will arrive at a more general summation formula
which we omit due to its cumbersome look.  A more general formula of this type will be given in our paper \cite[Theorem~2]{KPNew}.

\medskip

\textbf{Example~4.}  Take $r=3$, $\m=(1,1,2)$, $m=4$. Then formula (\ref{eq:KRPThKummerq-m-2}), in view of the remark after Corollary~\ref{cr:KRPTh1limitm-2},
takes the form
\begin{multline*}
e^{-x}{}_{4}F_{4}\left.\!\left(\!\begin{matrix}\bbbb,f_1+1,f_2+1,f_3+2\\\bbbb+2, f_1,f_2,f_3\end{matrix}\right\vert x\!\right)=
\sum\limits_{j=0}^{2}\frac{(-x)^j}{j!}{}_{5}F_{4}\!\left(\begin{matrix}-j,\bbbb,f_1+1,f_2+1,f_3+2\\\bbbb+2,f_1,f_2,f_3\end{matrix}\right)
\\
+(-x)^{3}B{}_{2}F_{2}\left.\!\!\left(\begin{matrix}1,\lambda+4\\\bbbb+5, \lambda+3\end{matrix}\right\vert-x\!\right),
\end{multline*}
where
$$
B=\frac{(\bbbb+4)(f_1-\bbbb)(f_2-\bbbb)(f_3-\bbbb)_2-\bbbb(f_1-\bbbb-1)(f_2-\bbbb-1)(f_3-\bbbb-1)_2}{(\bbbb+2)_{3}f_1f_2f_3(f_3+1)}
$$
and
$$
\lambda=1+\bbbb-\frac{\bbbb(f_1-\bbbb-1)(f_2-\bbbb-1)(f_3-\bbbb-1)}{(f_1-\bbbb)(f_2-\bbbb)(f_3-\bbbb+1)}.
$$
Theorem~\ref{th:Thomaelike2} takes the form ($q=2$):
$$
{}_{6}F_{5}\!\!\left(\begin{matrix}\aaaa,\bbbb,\dddd,f_1+1,f_2+1,f_3+2\\\bbbb+2, \eeee, f_1,f_2,f_3\end{matrix}\right)
=B_1-B_2{}_{4}F_{3}\!\left(\begin{matrix}1,\aaaa-\bbbb+1, \dddd+3, \gamma+4\\\aaaa+5, \eeee-\bbbb+1, \gamma+3\end{matrix}\!\right),
$$
where
\begin{multline*}
B_1=\frac{\Gamma(\eeee)\Gamma(\eeee-\bbbb-\dddd-2)}{\Gamma(\eeee-\dddd)\Gamma(\eeee-\bbbb-2)}\!\sum\limits_{j=0}^{2}\frac{(-2)_j(\aaaa-\bbbb-2)_j(\dddd)_j}{(\aaaa+2)_{j}(\eeee-\bbbb-2)_{j}j!}\times
\\
\sum\limits_{k=0}^{j}\frac{(-j)_{k}(\aaaa)_k(\bbbb)_k}{(-2)_k(\aaaa-\bbbb-2)_kk!}
{}_{4}F_{3}\!\!\left(\begin{matrix}-k,f_1+1,f_2+1,f_3+2\\f_1,f_2,f_3\end{matrix}\right)
{}_{3}F_{2}\!\!\left(\begin{matrix}k-4,k-j,-\bbbb-2\\k-2,\aaaa-\bbbb-2+k\end{matrix}\right),
\end{multline*}
$$
B_2=\frac{\Gamma(\eeee)\Gamma(\eeee-\bbbb-\dddd-2)(\bbbb-\aaaa)_{3}(\dddd)_{3}(\bbbb)_{3}[(\bbbb+3)(\f-\bbbb)_{\m}-\bbbb(\f-\bbbb-1)_{\m}]}{\Gamma(\eeee-\dddd)\Gamma(\eeee-\bbbb+1)(\aaaa+2)_{3}(\bbbb-\aaaa)_{3}(\f)_{\m}},
$$
$$
\gamma=\frac{\aaaa(\f-\aaaa-1)-(\aaaa+m)(\f-\aaaa-1+\m)}{\aaaa(\f-\aaaa-1)/(\aaaa-\bbbb+1)-(\f-\aaaa-1+\m)}-3
$$
and we utilized formula (\ref{eq:hatQm0negative}) when computing $B_1$ and the formula for $\hat{R}_{m-q-1}(1-m)$ from the proof of Corollary~\ref{cr:KRPTh2limitm-2} when computing $B_2$.
\bigskip
\bigskip

\textbf{Acknowledgements.} We express our gratitude to Professor Richard Paris for sharing his insights on the topic of this work.
This research has been supported by the Russian Science Foundation under the project 14-11-00022.

\end{document}